\documentclass[12pt,reqno,,usenames,dvipsnames]{amsart}
\usepackage[utf8]{inputenc}

\author[]{Nir Schwartz}
\address{Université Paris-Saclay, CNRS, Laboratoire de mathématiques d’Orsay, 91405, Orsay, France.}
\email{nir.schwartz@universite-paris-saclay.fr}
\usepackage{changes}
\definechangesauthor[name={Nir Schwartz}, color=blue]{N}
\usepackage{amssymb}
\usepackage{microtype}
\usepackage[shortlabels]{enumitem}
\usepackage{mathtools}
\usepackage{geometry}
\usepackage{soul}
\usepackage{mathrsfs}
\usepackage{stmaryrd}
\usepackage{setspace}
\usepackage{aliascnt}
\usepackage[scr=boondoxo]{mathalfa}
\usepackage{tikz}
	\usetikzlibrary{arrows, arrows.meta}
	\usetikzlibrary{decorations.pathmorphing}
	\tikzstyle{printersafe}=[decoration={snake,amplitude=0pt}]
\DeclareSymbolFont{yhlargesymbols}{OMX}{yhex}{m}{n}
\DeclareMathAccent{\wideparen}{\mathord}{yhlargesymbols}{"F3}

\makeatletter
\let\mymakefnmark\@makefnmark
\let\mythefnmark\@thefnmark

\newcommand{\restorefn}{\let\@makefnmark\mymakefnmark
	\let\mythfnmakr\@thefnmark}
\makeatother

\thanks{}
\allowdisplaybreaks
\theoremstyle{remark}
\usetikzlibrary{patterns}
\usepackage{scalerel,stackengine,wasysym}
\usepackage{scalerel}
\usepackage{stackengine}

\usepackage{enumitem}

\newcommand{\ca}[1]{\mathcal{#1}}

\newcommand{\bltr}{\mathcal{B}(\Ltr)}

\newcommand{\Ltr}{L^2(\RR)}

\newcommand{\hi}{\mathcal{H}_N}

\newcommand{\Omg}{\Omega}

\newcommand{\met}{\mathcal{M}_{N}}

\newcommand{\trs}{\pitchfork}
\newcommand{\TT}{\mathbb{T}}
\newcommand{\RR}{\mathbb{R}}
\newcommand{\NN}{\mathbb{N}}
\newcommand{\ZZ}{\mathbb{Z}}
\newcommand{\CC}{\mathbb{C}}
\newcommand{\QQ}{\mathbb{Q}}
\newcommand{\Op}{\mathrm{Op}}
\newcommand{\inp}[1]{\left\langle#1\right\rangle }
\newcommand{\hilb}{\mathcal{H}}
\newcommand{\al}{\alpha}
\newcommand{\lm}{\lambda}
\newcommand{\ha}{\frac 12}
\newcommand{\Tr}{\text{Tr}}
\newcommand{\qt}{\frac 14}
\newcommand{\cal}[1]{\mathcal{#1}}
\newcommand{\supp}{\text{supp}\ }
\newcommand{\bnd}{\partial}

\newcommand{\ka}{\kappa}
\newcommand{\be}{\beta}
\newcommand{\dl}{\delta}
\newcommand{\igrk}{\mathcal{Y}}
\newcommand{\F}{\mathcal{F}}
\newcommand{\hil}{\mathcal{H}}
\newcommand{\eps}{\varepsilon}
\usepackage{bbm}
\newcommand{\ind}{\mathbbm{1}}
\newcommand{\norm}[3]{\left\Vert #1 \right\Vert_{#2}^{#3}}
\newcommand{\qmatrix}[1]{\begin{pmatrix} #1 \end{pmatrix}}
\newcommand{\slt}[1]{\text{SL}_{2}({#1})}
\newcommand{\vphi}{\varphi}
\newcommand{\sub}{\subset}
\newcommand{\abs}[1]{\left|#1\right|}
\usepackage{hyperref}
\newcommand{\eq}[2]{\begin{equation}{#1}\label{#2}\end{equation}}
\newcommand{\8}{\infty}

\newcommand{\bigslant}[2]{{\raisebox{-.2em}{$#2$}\left\backslash\raisebox{.2em}{$#1$}\right.}}

\newcommand{\bt}{\epsilon}

\newcommand{\br}{\breve}

\usepackage[toc,page]{appendix}
\theoremstyle{definition}
\newtheorem*{thm*}{Theorem}
\newtheorem{thm}{Theorem}[section]
\newaliascnt{def}{thm}
\newtheorem{defn}[def]{Definition}
\newaliascnt{prop}{thm}
\newtheorem{prop}[prop]{Proposition}
\newaliascnt{rem}{thm}
\newtheorem{rem}[rem]{Remark}
\newaliascnt{example}{thm}
\newtheorem{exmple}[example]{Example}
\newaliascnt{clm}{thm}

\newaliascnt{cor}{thm}
\newtheorem{cor}[cor]{Corollary}
\newaliascnt{clc}{thm}

\newaliascnt{lem}{thm}
\newtheorem{lem}[lem]{Lemma}
\aliascntresetthe{lem}
\aliascntresetthe{clc}
\aliascntresetthe{cor}
\aliascntresetthe{clm}
\aliascntresetthe{example}
\aliascntresetthe{rem}
\aliascntresetthe{prop}
\aliascntresetthe{def}

\makeatletter
\renewcommand{\xRightarrow}[2][]{\ext@arrow 0359\Rightarrowfill@Statement of results{#1}{#2}}
\makeatother

\makeatletter
\newcommand{\dotDelta}{{\vphantom{\Delta}\mathpalette\d@tD@lta\relax}}
\newcommand{\d@tD@lta}[2]{%
	\ooalign{\hidewidth$\m@th#1\mkern-1mu\cdot$\hidewidth\cr$\m@th#1\Delta$\cr}%
}
\makeatother

\stackMath
\newcommand\reallywidehat[1]{%
	\savestack{\tmpbox}{\stretchto{%
			\scaleto{%
				\scalerel*[\widthof{\ensuremath{#1}}]{\kern-.6pt\bigwedge\kern-.6pt}%
				{\rule[-\textheight/2]{1ex}{\textheight}} 
			}{\textheight}%
		}{0.5ex}}%
	\stackon[1pt]{#1}{\tmpbox}%
}
\makeatletter
\def\resetMathstrut@{%
	\setbox\z@\hbox{%
		\mathchardef\@tempa\mathcode`\[\relax
		\def\@tempb##1"##2##3{\the\textfont"##3\char"}%
		\expandafter\@tempb\meaning\@tempa \relax
	}%
	\ht\Mathstrutbox@\ht\z@ \dp\Mathstrutbox@\dp\z@}
\makeatother
\begingroup
\catcode`(\active \xdef({\left\string(}
\catcode`)\active \xdef){\right\string)}
\endgroup
\mathcode`(="8000 \mathcode`)="8000

\usepackage{graphicx}
\usetikzlibrary{decorations.pathreplacing,decorations.markings}
\tikzset{
  on each segment/.style={
    decorate,
    decoration={
      show path construction,
      moveto code={},
      lineto code={
        \path [#1]
        (\tikzinputsegmentfirst) -- (\tikzinputsegmentlast);
      },
      curveto code={
        \path [#1] (\tikzinputsegmentfirst)
        .. controls
        (\tikzinputsegmentsupporta) and (\tikzinputsegmentsupportb)
        ..
        (\tikzinputsegmentlast);
      },
      closepath code={
        \path [#1]
        (\tikzinputsegmentfirst) -- (\tikzinputsegmentlast);
      },
    },
  },
  mid arrow/.style={postaction={decorate,decoration={
        markings,
        mark=at position .5 with {\arrow[#1,scale=2]{stealth}}
      }}},
}

\renewcommand{\psmallmatrix}[1]{\left(\begin{smallmatrix}#1\end{smallmatrix}\right)}
\usepackage{eso-pic}

\begin{document}

\title{The full delocalization of eigenstates for the quantized cat map}
\maketitle
\begin{abstract}
We consider the quantum cat map -- a toy model of a {quantized} chaotic system. We show that {its} eigenstates are fully delocalized { on $\TT^2$} in the semiclassical limit (or equivalently that each semiclassical measure is fully supported on $\TT^2$). We adapt the proof of a similar result proved for the eigenstates of $-\Delta_g$ on compact hyperbolic surfaces from \cite{DJ17}, relying on the fractal uncertainty principle in \cite{BD16}.
\end{abstract}

\section{Introduction}

Let us consider a classical chaotic dynamical system. One of the central objectives in quantum chaos is understanding the extent to which the quantum counterpart of the system behaves like it in the high frequency limit. We illustrate this with a fundamental widely-studied example: {l}et $(M,g)$ be a compact surface of constant negative curvature. The geodesic flow acts on $S^*M$ viewed as the phase space. It is well known that in this geometry the geodesic flow is Anosov, in fact it is uniformly hyperbolic, ergodic and mixing. Quantifying the flow results in the Laplace-Beltrami operator $-\Delta_g$ on $L^2(M)$, {which} generates the Schrödinger equation. {We expand} $L^2(M)$ as a direct sum of eigenspaces arising from a choice of an orthonormal eigenbasis $\{\vphi_j\}_j$. {The quantum-classical correspondence establishes a link between the classical flow and the quantum operator.}
This {connection} is given in terms of the Egorov theorem. The first results {on spatial distribution of high frequency eigenmodes} are due to {Shnirelman, Zelditch and Colin de Verdière (\cite{S74},\cite{Z87},\cite{CdV85})} stating that almost every $\vphi_j$ is {asymptotically} equidistributed on $M$ in the high frequency limit, {a property called quantum ergodicity}. A corollary of it is that there exists a set $J\sub\NN$ of density 1 such that for every $f\in C(M),$
\eq{
\int_M f(x) |\vphi_j(x)|^2dx \xrightarrow[j\in J]{j\to\8} \int_M f(x) dx.
}{cvg}

Rudnick and Sarnak have conjectured that the eigenmodes of $-\Delta_g$ satisfy the quantum unique ergodicity (QUE) property, meaning that $J=\NN$, i.e., \eqref{cvg} holds for the full sequence $\{\vphi_j\}$. Quantum ergodicity and QUE can be both rephrased by lifting the eigenmodes to the {distributions} $dW_j\in \mathcal{D}'(T^*M)$ {analogous to Wigner distributions on $T^*M$ {which are defined by 
\[
a\mapsto \inp{\Op_h(a)\vphi_j,\vphi_j},\qquad a\in C^\8(T^*M).
\]}} One can always extract a convergent subsequence $\{dW_{j_k}\}_{k\ge 1}$, converging {in the distributional sense} to a probability measure {on $S^*M$}. The weak-* limits of $\{dW_j\}_j$ are called semiclassical measures, denoted $\mu_{\text{sc}}$, and are invariant with respect to the geodesic flow. They represent the microlocalization of $\{\vphi_{j_k}\}_{k\ge 1}$ on the phase space. In these notations quantum ergodicity amounts to the existence of a subsequence {$J\sub \NN$}  of density 1 such that $dW_{j}\xrightharpoonup[j\in J]{j\to\8} \mu_{\text{Liouv}}$ on $S^*M$  and QUE means that  the full sequence weakly-converges to the Liouville measure. {We note that since the geodesic flow is Anosov, there exist many probability measures invariant with respect to it apart from Lebesgue, e.g., measures supported on a fractal invariant set or on closed geodesics.} 

While QUE remains open, one can wonder which of these myriad probability measures, being invariant under the geodesic flow, can be obtained {as weak-* limits of $\{dW_j\}_j$}. Given such a measure, one can also wonder where it is localized. There are {several} results concerning the constraints that semiclassical measures have to satisfy. 

Some of these constraints are expressed in terms of the Kolmogorov-Sinai entropy quantifying the information-theoretic complexity of a given invariant measure. This entropy gives information on the localization of the measure: the higher the entropy, the more delocalized the measure. For instance {the} Liouville measure is the measure  of maximal entropy $1$ whereas a delta measure on a closed geodesic is of minimal entropy $0$. Anantharaman proved in \cite{A08} that the Kolmogorov-Sinai entropy of every $\mu_{\text{sc}}$ is positive. Joint with Nonnenmacher they showed in \cite{AN2} that the entropy is in fact bounded from below by half the maximal entropy. These bounds mean that $\mu_{\text{sc}}$ cannot be "too" localized. Nevertheless there exist flow-invariant measures with high entropy which are supported on {invariant proper subsets}.{ A recent result due to Dyatlov and Jin ,\cite{DJ17}, states that any semiclassical measure is fully-supported on $S^*M$,{and moreover that} for every open $\emptyset\neq\Omg\sub S^*M$ there exists a constant $C_\Omg$ independent of the choice of $\mu_{sc}$ such that $\mu_{\text{sc}}(\Omg)>C_\Omg$. Their method of proof relies on the fractal uncertainty principle introduced in \cite{BD16} (see \autoref{propfup2}). In addition, combined with the unique continuation principle they deduce that for every open $\Omg\neq\emptyset$ there exists a constant $c_\Omg$ such that for every normalized eigenmode $\vphi_j$ \eq{\int_{\Omega}|\vphi_j|^2dx \ge c_\Omg\int_{M} |\vphi_j|^2dx=c_\Omg.}{bM}}
  
\begin{figure}[htp]
\centering
\includegraphics[scale=0.35]{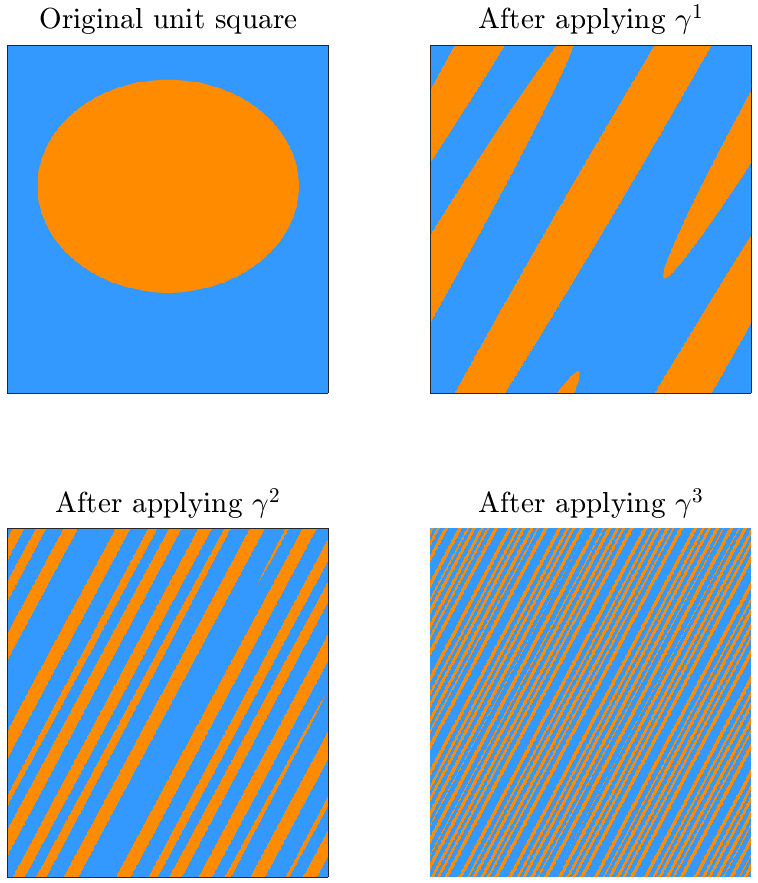}
\caption{Applying a hyperbolic toral automorphism $\gamma$ to a cluster of points (in orange) stretches them along the unstable branches, eventually filling {densely} the torus.}
\label{fig4}
\end{figure}

{In the present article} we study {the cat map,} a toy model of discrete-time dynamics on {the 2-dimensional torus} $\TT^2$, {having} analogous dynamical properties. This {cat map is an} Anosov diffeomorphism acting on $\TT^2$ and preserving its Liouville measure.
We consider the torus $\TT^2$ viewed as a symplectic manifold equipped with coordinates $(y,\eta)\in\RR^2/\ZZ^2$ as the phase space. The classical dynamics on it is given by a toral hyperbolic automorphism $\gamma\in\slt{\ZZ}$ which is an Anosov diffeomorphism (see \autoref{fig4}). A class of examples was popularized by Arnold in 1967 (cf. \cite{A67}) 
 and is named after him, Arnold's cat maps. Given a $\gamma\in\tilde\Gamma(2)<\slt{\ZZ}$ (where $\tilde\Gamma(2)$ is defined in \eqref{Gamma2}) Hannay and Berry  have quantized the map, i.e., {constructed} a family of unitary operators $\{\met(\gamma)\}_N$ of rank $N$, for $N\in\NN^*$ (see \autoref{qdyn} below and \cite{HB80}).{  This family of unitary operators is called the quantum cat map, associated with the symplectomorphism $\gamma$. Each $\met(\gamma)$ acts on an $N-$dimensional space of distributions $\hi\cong \CC^N$ recalled precisely in \autoref{qspace}.
 	
 	 We mention that each $\met(\gamma)$ can be viewed as an analogue of the propagator for the semiclassical Schrödinger equation, $e^{ih\Delta}$} {where $h$, being effective Planck's constant, is the semiclassical parameter}. {From a physical perspective, one can view $N$ as $N=\frac{1}{2\pi h}$,  i.e., the inverse of the effective Planck's constant. It means that while the semiclassical limit in the previous model corresponds to $h\to0$, in our model it corresponds to $N\to\8$.}

Each unitary operator $\met(\gamma)$ admits an orthonormal eigenbasis.  We note that high spectral degeneracy permits various choices of such basis. {We denote a normalized eigenvector of $\met(\gamma)$ by {$\vphi_N\in\CC^N$.}}  {The quantization procedure allows to associate to $\vphi_N$ a distribution $dW_N$ on $\TT^2$ analogous to {a} Wigner distribution}.{

Similarly one can construct semiclassical measures $\mu_{\text{sc}}$ as the weak-* limits of $\{dW_{N}\}$ . Each of them is a $\gamma-$invariant probability measures on $\TT^2$}. As for the setting of compact hyperbolic manifolds introduced above one can examine the {"delocalization"} of $\{\vphi_{N}\}$ or equivalently of {$\{dW_N\}$} in the semiclassical limit $N\to\8$.

Bouzouina and De Bièvre have proved in \cite{BD96} an analogue of quantum ergodicity for the quantum cat map, {i.e., there exists a set {$J\sub \NN$}  of density 1 such that $dW_{j}\xrightharpoonup[j\in J]{j\to\8} \mu_{\text{Liouv}}$}. However, quantum ergodicity does not prevent exceptional subsequences converging to a {$\gamma-$invariant} measure different from $\mu_{\text{Liouv}}$ on $\TT^2$. A work by Kurlberg and Rudnick (\cite{KR00}) {have constructed "joint" eigenbases composed} of mutual eigenvectors of $\met(\gamma)$ and "Hecke" operators on $\TT^2$  {(analogous to} those defined for compact arithmetic surfaces). They proved that these {eigenbases} satisfy quantum unique ergodicity.  In 2003, Faure, Nonnenmacher and De Bièvre (\cite{FNDB}) demonstrated that for the quantized cat map quantum unique ergodicity does not hold. For every periodic orbit $\mathcal{P}$ of $\gamma$, they found a weak-* limit of the form $\mu_{\text{sc}}=\frac12\mu_{\text{Leb}}+\frac 12\delta_{\mathcal{P}}$. Faure and Nonnenmacher later showed in \cite{FN4} that {for every $\mu_{\text{sc}}$ of the cat map} the weight of its atomic component is less or equal to the weight of its Liouville component, implying that the semiclassical measure cannot be entirely carried on a periodic orbit.
Anantharaman and Nonnenmacher have obtained a result on the entropy of semi-classical measures for another toy model on $\TT^2$, the Walsh-quantized Baker map (cf. \cite{An1}),namely that every semiclassical measure has at least half the maximal (i.e., Lebesgue) entropy. 
Brooks \cite{B10} has given a related characterization of semi-classical measures for the cat map. He interpolated between the constraint from \cite{FN4} and the analogue of the lower bound from \cite{An1} for the cat map. He proved that ergodic components of high entropy of the semiclassical measures carry at least as much weight as those of low entropy. 

Our paper is dedicated {to} proving a result analogous to \cite{DJ17} for the quantum cat map, that every semi-classical measure $\mu_{\text{sc}}$ {associated with $\gamma$} is fully supported on $\TT^2$ {meaning that all sequences of eigenstates $\{\vphi_{N}\}_N$ are fully delocalized {on $\TT^2$} for $N$ large enough, }
\begin{thm}\label{main}
Let $\gamma\in \tilde\Gamma(2)$ (where $\tilde\Gamma(2)$ is defined in \eqref{Gamma2}) {be a hyperbolic matrix quantized into the family $\{\met(\gamma)\}_N$. Let $\mu_{\text{sc}}$ be an associated semiclassical measure.} Then for every open $\emptyset\neq\Omg\sub\TT^2$ there exists a constant $c_\Omg>0$ independent of $\mu_{\text{sc}}$ such that $\mu_{\text{sc}}(\Omg)> c_\Omg$.
\end{thm}

We can deduce a corollary stated {in terms of $\vphi_N\in\CC^N$} {through the identification \eqref{iden} below}:
\begin{thm}[{The "delocalization" of the eigenfunctions}]\label{lacor}
Let $\gamma\in\tilde\Gamma(2)$ {and let $\vphi_{N}= (\vphi_{N,k})_{k=1}^{N}$ be a normalized eigenvector of $\met(\gamma)$. Fix $0\le \al_1 < \al_2\le 1$. Then there exists a constant $c_{\al_1,\al_2}>0$ and an index $N_0\in\NN^*$ such that for every $N>N_0$
	\eq{
		{ \frac 1N}\sum_{k\in \llbracket \al_1N,\al_2N\rrbracket} \abs{\vphi_{N,k}}^2 > c_{\al_1,\al_2}.
	}{sumcjkN}}
\end{thm}

{ Both theorems are deduced from \autoref{prop:blabla} (as we \hyperlink{pflacor}{prove in section 3}) .} 

We mention that one can {ask} similar questions for higher-dimensional analogues of the cat map hyperbolic dynamics, meaning hyperbolic maps $\gamma\in\text{Sp}(2d,\ZZ)$. Rivière (cf. \cite{R11}) has generalized the result of \cite{AN2} for these higher dimensional settings bounding explicitly the entropy of $\mu_{\text{sc}}$ from below in terms of the spectrum of $\gamma$. In addition Kelmer has obtained in \cite{K10} that if $\gamma$ belongs to a certain subgroup of $\text{Sp}(2d,\ZZ),$ one can {construct}  eigenvectors of the {corresponding quantum operators} $\met(\gamma)$  {for which the associated semiclassical measure is the Liouville} measure on an invariant Lagrangian subspace of $\TT^{2d}$, disproving QUE for $\gamma$ {as well as an analogue of \autoref{main} in high dimension}.

{ Recently Dyatlov and Jézéquel extended (cf. \cite{D21})  our result to cat maps of arbitrary dimension $\gamma\in\text{Sp}(2d,\ZZ)$. Their method of proof relies as well on the one dimensional fractal uncertainty principle (\autoref{propfup2} below).}

\subsection*{The structure of the paper}
{The next section is dedicated {to} recalling the quantum dynamics and the kinematics on $\TT^2$ and for {extending} the anisotropic pseudodifferential calculus introduced in \cite{DZ16} for the setting of the torus. In \autoref{quanmain} we introduce the main estimate in \autoref{prop:blabla}, from which we deduce our central results theorems \ref{main} and \ref{lacor}, as explained \hyperlink{pff}{below}. The proof of {\autoref{prop:blabla}} is done in several steps: We first set up a pseudo-differential partition of unity and then we refine it by propagating with respect to the quantum dynamics. Based on the refined partition we divide operators into two classes: a "controlled" class and an "uncontrolled" one. The proof of \autoref{prop:blabla} relies on estimating norms of operators from these classes: In \autoref{propigrk} we bound the norm of the controlled operators using {the anisotropic pseudodifferential calculus} introduced in \autoref{qdyn}. In the rest of the paper we consider the "uncontrolled" operators. In \autoref{fups} we prove that each symbol propagated for a long time has an almost fractal support. In \autoref{finalfups} we apply the fractal uncertainty principle introduced in \cite{BD16} estimating the norm of each operator. Combining the estimates for the norms of the controlled and uncontrolled operators implies \autoref{prop:blabla}.

\subsection*{Acknowledgments} This work was carried as part of the author's PhD research under the supervision of Stéphane Nonnnmacher, whose guidance made this project possible.{ The author would like to thank the anonymous referees for the constructive suggestions on the first versions of the paper, for Xiaolong Han and Gabriel Rivière for their discussions on entropy and capacity.} The author was supported by a MESRI grant.

\section{Quantum dynamics on the torus}\label{qdyn}
\subsection{Hyperbolic automorphisms}
{Let us consider $\RR^2$ equipped with the coordinates $(y,\eta)$. We call $y$ the "position coordinate" and $\eta$ the "momentum coordinate". We view the 2-dimensional torus as the quotient $\TT^2=\RR^2/\ZZ^2$.  The projection map between the two spaces will be denoted by 
	\eq{\pi(y,\eta) = (y,\eta) \pmod{\ZZ^2}\quad \pi:\RR^2\to\TT^2.}{pi}
	 We study the dynamics arising from a class of hyperbolic automorphisms, that is} $\gamma\in\slt{\ZZ}$ with eigenvalues $|\lm_u|>1,|\lm_s| < 1$ or equivalently $|\Tr(\gamma)|>2$.
 
 \begin{lem}\label{irr}
 	Let $v$ be an eigenvector of $\gamma$, then $v\in\RR^2\setminus\mathbb{Q}^2$.
 \end{lem}
 \begin{proof}
 	Let us first realize that $\lm_u\notin\mathbb{Q}$: Recall that $\text{Tr}(\gamma)^2-4$ is a square if and only if there exists $n\in\ZZ$ such that $\text{Tr}(\gamma)^2-n^2=4$, or equivalently if and only if $$\qmatrix{1&1\\1&-1}\binom{\text{Tr}(\gamma)}{n}\in\left\{\pm\binom 14,\pm\binom 41,\pm\binom 22\right\}.$$ {Multiplying both sides by $\qmatrix{1&1\\1&-1}^{-1}=\qmatrix{\ha&\ha\\\ha&-\ha}$}, in the first two cases $\text{Tr}(\gamma)\notin\ZZ$ and in the third $\text{Tr}(\gamma)=\pm 2$, contradicting the hyperbolicity of $\gamma$.
 	
 	Assume to the contrary that $v\in\mathbb{Q}^2$ is an eigenvector of $\gamma$, without loss of generality in the unstable direction. Then $\gamma v\in\mathbb{Q}^2$ whereas at least one component of $\lm_u v$ is irrational in contradiction. 
 \end{proof}
 
  {We deduce that the corresponding eigenspaces can be expressed by} \[W_u(\gamma)=\text{span}\left\{ \qmatrix{1\\m_u}\right\},\qquad W_s(\gamma)=\text{span}\left\{\qmatrix{1\\m_s}\right\}.\]
  {with $m_u,m_s\in\RR\setminus\QQ.$}
{Henceforth we refer to $W_u(\gamma)$ as the unstable space and to $W_s(\gamma)$ as the stable space.}
Strictly speaking, we are interested in hyperbolic automorphisms arising from the group
\eq{
\tilde\Gamma(2)=\left\{\qmatrix{a_{11}&a_{12}\\a_{21}&a_{22}}\in\slt{\ZZ}: a_{11}a_{12}\equiv a_{21}a_{22}\equiv 0\mod 2\right\}.
}{Gamma2}
{
We recall several elementary properties of elements in this group:
}

\begin{lem}
Let $\gamma\in\tilde\Gamma(2),$ then $\Tr (\gamma) \in 2\ZZ$.
\end{lem}

\begin{proof}
Let $\gamma=\psmallmatrix{a_{11}&a_{12}\\a_{21}&a_{22}}$. Assume that  $\Tr(\gamma)=a_{11}+a_{22}\in 2\ZZ+1$. Without loss of generality $a_{11}\in2\ZZ+1$ and $a_{22}\in2\ZZ$. Then $a_{11}a_{22}=\det(\gamma)+a_{12}a_{21} = 1+a_{12}a_{21}\in 2\ZZ$. As follows $a_{12},a_{21}\in2\ZZ+1$ and thus $a_{11}a_{12}\not\equiv 0\pmod 2$ contradicting $\gamma\in\tilde\Gamma(2).$
\end{proof}

\begin{cor}\label{mineig}
{For every hyperbolic matrix $\gamma\in\tilde\Gamma(2)$ from the previous lemma $\abs{\Tr(\gamma)}\ge 4$ and \[
\max_{\lm\in\text{Spec}(\gamma)}\abs{\lm}=\ha\cdot \max\left\{\abs{\Tr(\gamma)+\sqrt{\Tr(\gamma)^2-4}},\abs{\Tr(\gamma)-\sqrt{\Tr(\gamma)^2-4}}\right\}\ge 2+\sqrt 3.\] }
\end{cor}
{In addition, at least one of the coordinates of an eigenvector of $\gamma$ has to be irrational:}

{We end this subsection with a result on the minimality of the associated horocyclic flows: Let us fix an unstable eigenvector $v_u\in W_u$ and a stable eigenvector, $v_s\in W_s$. We define the unstable and the stable horocylic continuous flows associated to $\gamma$ by }
	\[
	\mathcal{H}_u(v,t) := v + tv_u\pmod{\ZZ^2} \qquad	\mathcal{H}_s(v,t) := v + tv_s\pmod{\ZZ^2} \qquad (v,t)\in \TT^2\times \RR 
	\]

We can verify by a direct calculation that both flows satisfy an intertwining relation with $\gamma,$
\eq{
	\gamma(\mathcal{H}_u(v,t)) = \mathcal{H}_u(\gamma.v,\lm_u t)\qquad 	\gamma(\mathcal{H}_s(v,t)) = \mathcal{H}_s(\gamma.v,\lm_s t)
}{Hin}
From \autoref{irr} we deduce that both flows are minimal,
\begin{lem}\label{dens}
	Every orbit $\mathcal{O}_{u}$ of $\mathcal{H}_u$ and every orbit $\mathcal{O}_{s}$ of $\mathcal{H}_s$ is dense in $\TT^2$.
\end{lem}
For sake of completeness we recall the proof.
\begin{proof}
	The proof for orbits of $\mathcal{H}_s$ is completely analogous thus we only prove the lemma for the orbits of $\mathcal{H}_u$. Let $(y_0,\eta_0)\in\TT^2$ and fix an unstable orbit $\mathcal{O}_1=\{\binom{y_1}{\eta_1} + tv_u \pmod \ZZ^2\}$ where $v_u=\binom{1}{m_u}$ as before. Let us consider the set
	\[
	\mathcal{T} := \{t\in\RR: t - (y_0-y_1) \in \ZZ \}.
	\] 
	and the corresponding sequence of points on the orbit $\{\binom{y_0}{m_un + \eta_1+m_u(y_0-y_1)}\pmod{\ZZ^2} \}_{n\in\ZZ}\sub\mathcal{O}_1\cap\{y\equiv y_0\pmod{\ZZ}\}$. Then, $\{m_un\pmod\ZZ\}_{n\in\ZZ}$ is a dense sequence on $[0,1]$ as $m_u\in\RR\setminus\QQ$ hence we can extract a sub-sequence of times ${t_{n_k}}$ for which $\mathcal{H}_u\left(\binom{y_1}{\eta_1},t_{n_k}\right)\xrightarrow{k\to\8} \binom{y_0}{\eta_0}$.
	\end{proof}

{We recall below that each $\gamma\in\tilde{\Gamma}(2)$ is quantized by a family of unitary operators $\{\met(\gamma)\}$ satisfying {an exact} Egorov property \eqref{egorovGamma}.} 
\subsection{Anisotropic calculi associated to Lagrangian foliations}
{Let us recall the semiclassical parameter $h\in(0,1)$ and recall {the class of symbols} {
 \begin{align*}S(\RR^2)=\{a=a(h)\in C^\8(\RR^2): &|\bnd_y^k\bnd_\eta^m a(y,\eta;h)| < C_{k,m}(a),\\ &\qquad\qquad C_{k,m}(a)\in\RR \text{ depending only on }k,m,a \}.\end{align*}}
We would like to propagate smooth symbols $a\in C^\8(\RR^2)$ with respect to the dynamics of $\gamma$, up to {large semiclassical time, namely twice the Ehrenfest time, that is denoting the semiclassical parameter by $h>0$ we propagate up to time} $2\left\lfloor\frac{\log \frac 1h}{\log \lm_u}\right\rfloor$. That results in highly oscillatory symbols lying outside $S(\RR^2)$, hence one has to define classes of symbols smooth along a 1-dimensional subspace of $\RR^2$ and oscillating along another.}

In our settings we consider only linear foliations, that is those having as leaves $L_{(y_0,\eta_0)}$ lines with a fixed incline passing through $(y_0,\eta_0)$. By a direct calculation one dimensional foliations are {always} Lagrangian.

\begin{rem}\label{Lsu} Each hyperbolic toral automorphism $\gamma$ induces two one-dimensional foliations, an unstable $L^{u}$ and a stable $L^{s}$, both having linear leaves.
Define the global unstable leaf at $(y_0,\eta_0)\in\RR^2$ by $W^{\RR^2}_{u,(y_0,\eta_0)}=\{(y,\eta)\in\RR^2: \eta=m_uy+\eta_0-m_uy_0 \}$ and similarly the global stable leaf by $W^{\RR^2}_{s,(y_0,\eta_0)}=\{(y,\eta)\in\RR^2: \eta=m_sy+\eta_0-m_sy_0 \}$. These yield the unstable and stable foliations of a hyperbolic $\gamma\in\slt{\RR}$, namely the foliations $L^u,L^s$ with leaves which are explicitly given by $W^{\RR^2}_{u,(y_0,\eta_0)}$ and $W^{\RR^2}_{s,(y_0,\eta_0)}$ correspondingly.
\end{rem}
Let $L$ be a linear foliation and let $v\in\RR^2$ be the vector associated to it. Fix a $v^\trs\in\RR^2$ transversal to $v$. The foliation associated to $v^\trs$ is denoted by $L_{v^\trs}$ and is called a $L-$transversal foliation.
{We remind that a symbol is a $h$-dependant function $a(h)=a(y,\eta;h)\in C^\8(\RR^2)$. } We introduce the notation
\[
a(h) =  O(h^{\zeta-}) \text{ if } a(h) = O(h^{\zeta-\eps}) \text{for all $\eps>0$}.
\]
Let us introduce two such of symbols independent of the choice of $v^\trs$,
\begin{defn}
Let $L$ be a linear foliation, $L_{v^\trs}$ a $L-$transverse foliation and $V_L, V_{L_{v^\trs}}$ vector fields acting by a directional derivative with respect to $v$ and $v^\trs$ correspondingly. {A symbol} $a\in S_{L,\rho+}\left(\RR^2\right)$ if {$a:\RR^2\to \RR$} is smooth and $k,m\in\NN$ for every $h\in(0,1)$ \[\sup_{(y,\eta)\in U}\left|V_L^m V_{L_{v^\trs}}^ka\left(y,\eta;h\right)\right|=O(h^{-\rho k-}) .\]

\end{defn}

This class is larger than the one introduced in \cite[Definition 3.2]{DZ16} as we ask the bound above to hold for every $\eps$ whereas in \cite{DZ16} the authors fix $\eps=0$.

\begin{rem}\label{remsl}
	We note that these classes are independant of the choice of $L_{v^\trs}$: given two different $L-$transverse foliaitons $L_{u^{\trs}},L_{v^\trs}$, we can write $V_{L_{u^\trs}} = C V_L + C' V_{L_{v^\trs}}$ and using the triangle inequality\[
	\abs{V_L^mV^k_{L_{u^\trs}}a(y,\eta;h)} \le \sum_{j+j'=k} C_{j,j',u^\trs,v^\trs}	\abs{V_L^{m+j}V^{j'}_{L_{v^\trs}}a(y,\eta;h)}  = O(h^{-\rho k-}).
	\]
	A similar argument shows that keeping $L_{v^\trs}$ while changing $L$, does not preserve the symbol class. As follows, we realize symbols in these anisotropic classes have sharp oscillations along any direction transerve to $L$ and controlled oscillations along $L$.
\end{rem}

For future reference let us recall the isotropic classes \begin{align*}
S_\rho(\RR^2)  &=\{a\in C^\8(\RR^2):\forall\ka\in\NN^2,\quad |\bnd^{\kappa}a|= O(h^{-\rho |\ka|})\}.\\
S_{\rho+}(\RR^2) &= \{a\in C^\8(\RR^2):\forall\ka\in\NN^2,\quad |\bnd^{\kappa}a|= O(h^{-\rho |\ka|-})\}\end{align*}
\begin{exmple}\label{exR2}
Let $\gamma\in\tilde\Gamma(2)$ and denote $v_{u}=\binom{1}{m_u}, v_s=-\binom{1}{m_s}$ the vectors corresponding to its unstable and stable foliations (see \autoref{Lsu}). {Since} $v_s$ is transversal to $v_u$ we define the symbol classes
 \begin{align} \label{sym_class}
 S_{L^u,\rho+}(\RR^2) &= \left\{a\in C^\8(\RR^2):\forall k,m\in\NN,\quad \sup_{(y,\eta)\in \RR^2}\left|V_u^m V_s^k a\left(y,\eta;h\right)\right| = O(h^{-\rho k-})\right\}, \\ S_{L^s,\rho+}(\RR^2) &= \left\{a\in C^\8(\RR^2):\forall k,m\in\NN,\quad \sup_{(y,\eta)\in \RR^2}\left|V_s^m  V_u^k a\left(y,\eta;h\right)\right| =  O(h^{-\rho k-}) \right\}\nonumber.
 \end{align}

\end{exmple}

\begin{defn}\label{VTT}
Let us consider analogous symbol classes on $\TT^2$. Every linear foliation $L$ can be projected to a foliation $L_{\TT^2}$ on $\TT^2$. {Equivalently} $V_{L_{\TT^2}}$ is the vector field on $\TT^2$ acting on $\ZZ^2$-periodic functions $a\in C^\8(\RR^2)$ by
\[
V_{L_{\TT^2}}a(y,\eta;h) := {V_L a(y,\eta;h)}.
\] 
Let $\gamma\in\tilde \Gamma(2)$ acting on $\TT^2$ and let $v_u,v_s$ as above. We denote the associated linear foliations on $\TT^2$ as $L^u_{\TT^2}$ and $L^s_{\TT^2}$ {}. We denote by $C^\8(\TT^2)$ the class of {smooth real-valued functions on $\TT^2$}. The symbol classes on the torus are defined as,
 \begin{align} \label{sym_class}
 S_{L^u,\rho+}(\TT^2) &= \left\{a\in C^\8(\TT^2):\forall k,m\in\NN,h\in(0,1)\quad \sup_{(y,\eta)\in \TT^2}\left|V_{L_{\TT^2}^u}^m V_{{L_{\TT^2}^s}}^k a\left(y,\eta;h\right)\right|= O(h^{-\rho k-})\right\}, \\ S_{L^s,\rho+}(\TT^2) &= \left\{a\in C^\8(\TT^2):\forall k,m\in\NN,h\in(0,1),\quad \sup_{(y,\eta)\in \TT^2}\left|V_{L_{\TT^2}^s}^m  V_{L_{\TT^2}^u}^k a\left(y,\eta;h\right)\right| = O(h^{-\rho k-})\right\},\nonumber
 \end{align}
with $V_{L_{\TT^2}^u}$ and $V_{L_{\TT^2}^s}$ being the vector fields associated to the linear foliations $L^u_{\TT^2}$ and $L^s_{\TT^2}$.
\end{defn}

The next class allows stronger oscillations along $L$,
\begin{defn}\label{rhorho'}
	Let $L$ be a linear foliation, $V_L, V_{L_{v^\trs}}$ vector fields acting by a directional derivative with respect to $v$ and $v^\trs$ correspondingly and let $\rho,\rho'\in \left[0,1\right)$ satisfying 
	\[
	0\le \rho' < \rho\qquad \rho+\rho'<1.
	\]
	{A symbol }$a\in S_{L,\rho,\rho'}(\RR^2)$ if $a:\RR^2\to\RR$ is smooth and for every $k,m\in\NN$ there exists an $h$-independent constant $C_{k,m}$ such that for every $h\in(0,1)$ \[\sup_{(y,\eta)\in \RR^2}\left|V_L^ m V^k_{L_{v^\trs}}a\left(y,\eta;h\right)\right|\le C_{k,m} h^{-\rho k-\rho'm}.\]
\end{defn}

\begin{rem}
	As in \autoref{remsl} the class characterizes the direction $L$ rather than $L_{v^\trs}$. In fact, symbols in the class can have sharper oscillations along the transversal direction $L_{v^\trs}$.
\end{rem}

\begin{exmple} For $L=L^{\eta}$ the vertical foliation in $\RR^2$ fixing as the transversal foliation $L_{v^\trs}=L^y$ we obtain a model {$\RR^2$-}calculus slightly broader compared to the one introduced in [Appendix A.2, \cite{DJ17}] which allowed only compact symbols:
\end{exmple} 

\begin{lem}\label{lemsym} Let $0\le 
	\rho,\rho' <1$.{ The class of symbols $S_{L^\eta,\rho,\rho'}\left(\RR^2\right)$
		consists of the functions $a\in C^{\8}\left(\RR^2\right)$ such that for every $k,m\in\NN$ there exists a constant $C_{k,m}>0$ satisfying for every $h\in(0,1)$ the derivative bounds \eq{\sup_{\left(y,\eta\right)\in\RR^{2}}\left|\bnd_{y}^{k}\bnd_{\eta}^{^{m}}a\left(y,\eta;h\right)\right|\le C_{k,m}h^{-\rho k-\rho'm}.}{a_hor} }
\end{lem}
\begin{proof}
	{ Any vector associated to $L^\eta$ is of the form $\binom{0}{j}$ for $j\in\RR^*$ and any vector associated to $L_{v^\trs}$ is of the form $v^\trs=\binom{j'}{0}$ with $j'\in\RR^*$. Therefore for every $k,m\in\NN$ there is a constant $C_{k,m}$ satisfying
		\[
		\sup_{(y,\eta)\in \RR^2}\abs{\bnd_\eta^m \bnd_y^k a(y,\eta;h)} = j^{-m}(j')^{-k}\sup_{(y,\eta)\in\RR^2}\abs{V_{L^\eta}^m V_{L^y}^k a(y,\eta;h)}\le C_{k,m} h^{-\rho k -\rho'm}.
		\] }
\end{proof}
{ We note that for every foliation $L$
\eq{
	S_{L,\rho+}(\RR^2) = \bigcap_{\eps>0} S_{L,\rho+\eps,\eps}(\RR^2)
}{}
hence \autoref{lemsym} implies that

\begin{align*}
S_{L^\eta,\rho+}(\RR^2) = \{a\in C^\8(\RR^2): \text{For every } k,m\in\NN  &\text{ exists a constant } C_{k,m,\eps}\text{ such that }\\ &\sup_{(y,\eta)\in \RR^2}\abs{\bnd_y^k\bnd_\eta^m a(y,\eta;h)} = O(h^{-\rho k -})\}.
\end{align*}
}

There is an algebraic identity which reduces proofs for {these classes of symbols} to the model calculus. Recall that
\[
m_u = \frac{a_{22}-a_{11}+\sqrt{\Tr^2(\gamma)-4}}{2a_{12}} \text{ or } m_u = \frac{a_{22}-a_{11}-\sqrt{\Tr^2(\gamma)-4}}{2a_{12}}.
\] From the the hyperbolicity of $\gamma,\ m_u\neq m_s$. We assume henceforth that $m_u > m_s$ (as in \autoref{excat} below).
Define $\iota\in\slt{\RR}$ by 
\eq{
\iota =\frac{1}{ \sqrt{m_u-m_s}} \begin{pmatrix} 1&-1\\ m_u&-m_s\end{pmatrix},\qquad \iota^{-1} = \frac{1}{\sqrt{m_u-m_s}}\begin{pmatrix}
-m_s &1\\-m_u& 1
\end{pmatrix}.
}{iota}

Note that the columns of $\iota$ are obtained by rescaling $v_u,v_s$, ensuring $\iota\in\slt{\RR}$. This automorphism maps the standard orthonormal basis of $\RR^2$ to an eigenbasis of $\gamma$: \[\iota\qmatrix{1\\ 0} =\frac{1}{\sqrt{m_u-m_s}} v_u,\qquad \iota  \qmatrix{0\\1} = \frac{1}{\sqrt{m_u-m_s}}v_s .\]
{We note that $\iota$ maps every point $(y,\eta)\in U\sub \RR^2$ with $U$ being open to  $\binom x\xi:=\iota\binom y\eta\in \iota U$, hence $\iota$ can be viewed as a change from the {$(y,\eta)$} coordinate frame to the {$(x,\xi)$} stable/unstable coordinate frame.} 
{ Henceforth for every symbol $a\in C^\8(\RR^2)$ we define the twisted symbol $\widetilde{a}$ by $\widetilde{a} := a\circ\iota$.} 

\begin{lem}\label{lemLLxi}
Let $\gamma\in\tilde\Gamma(2)$ and $a\in S(\RR^2)$ then for $h\in(0,1)$ and uniformly in $t\in\left[0, \rho\frac{\log\frac 1h}{\log |\lm _u|}\right]$
\begin{enumerate}
\item $\widetilde{a\circ\gamma^t}\in S_{L^\eta,\rho,0}(\RR^2)$.
\item $\widetilde{a\circ\gamma^{-t}}\in S_{L^y,\rho,0}(\RR^2)$.
\item $a\circ\gamma^{t}\in S_{L^s,\rho,0}(\RR^2)$.
\item $a\circ\gamma^{-t}\in S_{L^u,\rho,0}(\RR^2)$.
\end{enumerate}
\end{lem}
\begin{proof}
We prove only (i) and (iii) as the others follow from replacing $\gamma$ by $\gamma^{-1}$ exchanging between the contracting and expanding eigendirections.
\begin{enumerate}
\item Suppose $a\in S(\RR^2)$ so does $a\circ\iota$ since {\eq{|\bnd_y^k\bnd_\eta^m \widetilde{a}(y,\eta)| \le B(k,m)\norm{\iota}{}{k+m} \max_{0\le k'+m'\le k+ m}|\bnd_x^{k'}\bnd_\xi^{m'} a(x,\xi;h)|\le C'_{k,m}(a)}{aiota}} with $B(k,m)$ being a combinatorial constant depending on $k,m$. We note that $\iota^{-1}$ diagonalises $\gamma$ hence
\eq{
\iota^{-1}\gamma^t\iota = D^t = \qmatrix{\lm_u^t&0\\0&\lm_u^{-t}}.
}{diag}
Multiplying by $\iota$ from the left and using \eqref{aiota} we conclude that for $0<t\le\left\lfloor\frac{\log\frac 1h}{\log |\lm_u|}\right\rfloor$, {
\begin{align*}
 \sup_{(y,\eta)\in\RR^2}\abs{\bnd_{y}^k\bnd_\eta^m (\widetilde{a\circ\gamma^t})(y,\eta;h)} &= \sup_{(y,\eta)\in\RR^2}\abs{\bnd_{y}^k\bnd_{\eta}^m (\widetilde{a} (\lm_u^t y,\lm_u^{-t}\eta;h))} \\
&= |\lm_u|^{(k-m)t} \sup_{(y,\eta)\in\RR^2}\abs{\bnd_{y'}^k\bnd_{\eta'}^m (\widetilde{a})(y',\eta';h)\mid_{(y',\eta')=(\lm_u^ty,\lm_u^{-t}\eta)}} \\&\le C_{k,m}'(a)|\lm_u|^{(k-m)t} \\
&\le C_{k,m}'(a)|\lm_u|^{k t}\\
&< C_{k,m}'(a) |\lm_u|^{k\rho\frac{\log\frac 1h}{\log |\lm_u|}} \\
& = C_{k,m}'(a)h^{-k\rho}.
\end{align*}}
Thus,
\[
\sup_{(y,\eta)\in\RR^2}\abs{\bnd_{x}^k\bnd_\xi^m \widetilde{a\circ\gamma^t}(y,\eta;h)} \le C_{k,m}(a) h^{-\rho k}.
\]
\setcounter{enumi}{2}
\item {We recall that for every symbol $b\in C^\8(\TT^2),\quad$ $(V_ub)\circ\iota=\bnd_y(b\circ \iota)$ and similarly $(V_sb)\circ\iota=\bnd_\eta(b\circ\iota)$ and therefore taking $b=a\circ \gamma^t$ we obtain from (i)
\[
\sup_{(y,\eta)\in\RR^2}\abs{V_u^kV_s^m(a\circ\gamma^t)(y,\eta)} = \sup_{(x,\xi)\in\RR^2}\abs{\bnd_y^k\bnd_\eta^m(\widetilde{a\circ\gamma^t})(x,\xi)} \le C'_{k,m}(a)h^{-\rho k}.
\]

}
\end{enumerate}
\end{proof}
\subsection{The quantum settings}
For every automorphism $\gamma\in\slt{\RR}$ we recall $\mathcal{M}_h(\gamma):\mathcal{S}(\RR)\to\mathcal{S}(\RR)$, the projective metaplectic representation of $\slt{\RR}$ which corresponds to $\gamma$ (cf. \cite{Fol16}). Since Schwartz space is dense inside $\Ltr$, one can extend $\mathcal{M}_h(\gamma)$ to a unitary operator $\mathcal{M}_h(\gamma):\Ltr\to\Ltr$ (or by a similar argument to $\mathcal{M}_h(\gamma):\mathcal{S}'(\RR)\to\mathcal{S}'(\RR)$).
{
We quantize symbols from all classes mentioned above by the Weyl quantization {(cf. subsection 4.1.1 in \cite{Z12})} $\Op_h^w(\cdot)=\Op_h(\cdot)$. This choice of quantization is motivated by the exact intertwining relations \eqref{egorov_full} and \eqref{egorov} which hold in our settings as we mention below. We recall that for every symbol $a\in S(\RR^2)$ the operators $\Op_h(a)$ map $\mathscr{S'}(\RR^2)$ to itself and are bounded on $\Ltr$ uniformly in $h$. The same statement holds for $a\in S_{L,\rho+}(\RR^2)$ as we show in \ref{SSa2} below.}
 By Stone-Weierstrass theorem the Weyl quantization of smooth bounded symbols satisfies an exact Egorov relation (see subsections  1.5 and 4.2 in \cite{Fol16}): {for every fixed $h$ and bounded $a\in S(\RR^2)$ denote by $\mathcal{M}_h$ the full unitary metaplectic representation of { $\text{Mp}_2(\RR)$} (playing the role of a double cover of $\slt{\RR}$). Identifying each $\gamma\in\tilde\Gamma(2)$ with an antecedent in $\text{Mp}_2(\RR)$,
\eq{
\Op_h(a)=\mathcal{M}_h(\gamma)\Op_h(a\circ\gamma)\mathcal{M}_h(\gamma)^*
}{egorovGamma}
and 
\eq{
\Op_h(a) =\mathcal{M}_h(\iota)\Op_h(\widetilde{a})\ca{M}_h(\iota)^*,
}{egorov_full}
independently of the size of the derivatives of $a$.} This intertwining relation between $\Op_h(a)$ and $\Op_h(\widetilde{a})$ plays a major rule in studying the properties of the calculus. 

\begin{lem} \label{slr+}
Let $a\in S(\RR^2)$ and $\mathscr{a}\in S_{L^s,\rho+}(\RR^2)$ {be real valued symbols}. Fix a choice of $\iota\in\slt{\RR}$ mapping the $(y,\eta)$ coordinate frame to the stable/unstable coordinate frame and by an abuse of notation $\tilde{\mathscr{a}}=\mathscr{a}\circ\iota$. Then for every $0 < t <\rho\frac{\log\frac 1h}{\log|\lm_u|}$
\begin{enumerate}[ref =\autoref{slr+}.\roman*]
\item\label{Ssa} Self-adjointness: $(\Op_h(a))^*=\Op_h(a)$ {on $L^2(\RR)$}.
\item \label{slr+:moyal} The Moyal product of symbols: if $d\in S_{L^s,\rho+}(\RR^2)$ then {$\Op_{h}(\mathscr a)\Op_{h}(d)=\Op_{h}\left(c\right)$
with $c=\mathscr{a}\cdot d+O_{S_{L^s,\rho+}(\RR^2)}(h^{1-\rho-}).$ In fact,
recalling \autoref{lemLLxi}.(iii), $\Op_{h}(a\circ\gamma^t)\Op_{h}(d)=\Op_{h}\left(c\right)$
with $c=(a\circ\gamma^t)\cdot d+O_{S_{L^s,\rho+}(\RR^2)}(h^{1-\rho-}).$}
\item\label{SSa2} The operator $\Op_h(\mathscr{a})$ is bounded on $\Ltr$ uniformly in $h$.
\item\label{slr+:garding} G\r{a}rding inequality: if $ d\ge0$ and $d\in S_{L^s,\rho,\rho'}(\RR^2)$ then for some constant $C$ independent of $h$, for every $\vphi\in L^2(\RR)$ $$\left\langle \Op_{h}(d)\vphi,\vphi\right\rangle\ge-Ch^{1-\rho-\rho'}\Vert\vphi\Vert_{L^{2}(\RR)}^{2}.$$
\end{enumerate}
\end{lem}

\begin{proof}
\begin{enumerate}
\item[] 

\item The Weyl quantization satisfies $(\Op_h(a))^*=\Op_h(a^*)$ and $a$ is a real symbol.
\item We introduce the unitary rescaling operator
\eq{
	T_{\rho}:L^2(\RR)\to L^2(\RR),\qquad (T_{\rho} f)(x')= h^{\frac{\rho}{4}} f(h^{\frac{\rho}{2}}x')
}{rescop}
{We note that both for symbols $a\in S(\RR^2)$ and for symbols $\widetilde{\mathscr{a}}\in S_{L^\xi,\rho{+}}(\RR^2)$ conjugation by $T_\rho$ leads to
\eq{T_{\rho}\Op_h(a)T_{\rho}^{-1}=\Op_h(a^\flat)}{resym} with $ a^\flat(x,\xi;h)= a(h^{\frac{\rho}{2}}x,h^{-\frac{\rho}{2}}\xi)\in S_{\frac{\rho}{2}}(\RR^2)$} (cf. subsection 4.4 in \cite{Z12}) {and $\widetilde{\mathscr{a}}^\flat\in S_{\frac\rho 2+}(\RR^2)$ since for every $n=(n_1,n_2)\in\NN^2$
\[
\sup_{(x,\xi)\in\RR^2}\abs{\bnd^{n} \mathscr{a}^\flat} = h^{\frac\rho 2(n_1-n_2)}\sup_{(x,\xi)\in\RR^2}\abs{\bnd^n \mathscr{a}} \le C_n h^{-\frac{\rho}{2}\abs{n}-\eps}.
\]
}Denoting {$b=\widetilde{\mathscr a}\in S_{L^\xi,\rho+}$} and applying theorem 4.18 from \cite{Z12} implies that 

\[
\Op_h({b}^\flat)\Op_h((\widetilde{d})^\flat) = \Op_h((\widetilde{c})^\flat)
\]

with $(\widetilde{c})^\flat={b}^\flat\cdot (\widetilde{d})^\flat+{O_{S_{\frac \rho 2}+}}(h^{1-\rho-})$.
{Conjugating $\Op_h(\cdot)$ by $\mathcal{M}_h(\iota)$ maps symbols from $S_{L^s,\rho+}$ to $S_{L^\xi,\rho+}$.} Then, by the exact Egorov relation \eqref{egorov_full} applied with $\iota^{-1}$, {
\begin{align*}
\Op_h(\mathscr a)\Op_h(d)& = \mathcal{M}_h(\iota)\Op_h(b)\mathcal{M}_h(\iota)^*\mathcal{M}_h(\iota)\Op_h(\widetilde{d})\mathcal{M}_h(\iota)^* \\
&= \mathcal{M}_h(\iota)\Op_h(b)\Op_h(\widetilde{d})\mathcal{M}_h(\iota)^* \\
&= \mathcal{M}_h(\iota)T_{\rho}^{-1}\Op_h(b^\flat)\Op_h((\widetilde{d})^\flat)T_{\rho}\mathcal{M}_h(\iota)^* \\
&= \mathcal{M}_h(\iota)\Op_h(\widetilde{c})\mathcal{M}_h(\iota)^* \\
&= \Op_h(c).
\end{align*} 
}
\item Note that $\Op_h(\mathscr{a})=T_\rho^{-1}\Op_h(\mathscr{a}^\flat)T_\rho$. Since $\mathscr{a}^\flat\in S_{\frac\rho2+}(\RR^2)$ the operator $\Op_h(\mathscr{a}^\flat)$ and thus $\Op_h(\mathscr{a})$ are both bounded on $\Ltr$ uniformly in $h$.
\item First, from a direct calculation analogous to the one appearing in \autoref{lemLLxi} we deduce that $\tilde d\in S_{L^{\xi},\rho,\rho'}(\RR^2)$. Then, we recall the rescaling from \cite[Lemma A.2]{DJ17}, 
\[
h^\sharp:=h^{1-\rho-\rho'}\quad (\tilde d)^\sharp(x,\xi;h) := \tilde d(h^\rho x,h^{\rho'}\xi;h)\quad u^\sharp(x):=h^{\frac \rho 2}u(h^\rho x)
\]
for which using the exact intertwining relation in \eqref{egorov_full} \[
\inp{\Op_h(\tilde d)\mathcal{M}_h(\iota)\vphi,\mathcal{M}_h(\iota)\vphi} = \inp{\Op_{h^\sharp}((\tilde d)^\sharp)(\mathcal{M}_h(\iota)\vphi)^\sharp,(\mathcal{M}_h(\iota)\vphi)^\sharp} =\inp{\Op_h(d)\vphi,\vphi}
\]
and adapting \cite[Lemma A.4]{DJ17} we can also verify that $(\tilde d)^\sharp\in S(\RR^2)$. Then applying \cite[Theorem 4.32]{Z12} and using the unitarity of $\mathcal{M}_h(\iota)$ we deduce the inequality.
\end{enumerate}
\end{proof}

{The same argument used to prove \ref{slr+:moyal} can be applied to give a Moyal product in $S_{L,\rho,\rho'}(\RR^2)$ as well,
\begin{lem}[Moyal product of symbols in $S_{L,\rho,\rho'}(\RR^2)$]\label{moyalrhorhop}
	Fix a linear foliation $L$ and let $\mathscr{a},d\in S_{L,\rho,\rho'}(\RR^2)$ be real valued symbols, then 
	\(
		\Op_h(\mathscr{a})\Op_h(d) = \Op_h(c)
	\)
	with $c=\mathscr{a}\cdot d + O_{S_{L,\rho,\rho'}(\RR^2)}(h^{1-\rho-\rho'})$.
\end{lem}
}

We deduce a bound on $\norm{\Op_h(a)}{\Ltr}{}$ for {a bounded $a\in S_{L,\rho+}(\RR^2)$ {with $L$ being any linear foliation},
\begin{lem}\label{bdopa}
Suppose $a\in S_{L,\rho+}(\RR^2)$. If $\sup |a| \le 1$, for every $\eps>0$,
\[
\norm{\Op_h(a)}{}{} \le 1 + Ch^{1-\rho-\eps}.
\]
\end{lem}
}
\begin{proof}
In order to simplify the proof assume $a$ is a real symbol. We note that from \ref{slr+:garding}, for every $\vphi\in L^2(\RR)$
\[
\inp{\vphi,\vphi} - \inp{\Op_h(|a|^2)\vphi,\vphi} \ge -Ch^{1-\rho-\eps}\norm{\vphi}{}{2},
\]
{and as $\Op_h(a)\Op_h(a)=\Op_h({\abs{a}^2})+O_{\bltr}(h^{1-\rho-\eps})$}
\[
\norm{\Op_h(a)\vphi}{}{2}\le \norm{\vphi}{}{2}+Ch^{1-\rho-\eps}\norm{\vphi}{}{2},
\]
implying the lemma.
\end{proof}

\subsection{Example: Degli—Esposti's cat map}\label{excat}
We will present an example of hyperbolic automorphisms by Degli-Esposti's cat map (cf. \cite{DG03}), \[\gamma_{DE}=\qmatrix{2&1\\3&2}\in \tilde\Gamma(2)\] acting on $\RR^2$ (or rather on $\TT^2=\RR^2/\ZZ^2$). 
Its eigenvalues are $\lm_u=2+\sqrt 3,\ \lm_s=2-\sqrt 3$.  
For Degli—Esposti's cat map the corresponding eigenspaces are given explicitly by parallel transport of the vectors 
\[
v_u(\gamma)=\binom{1}{\sqrt 3}\qquad v_s(\gamma)=\binom{-1}{\sqrt 3}.\] 

A possible choice of $\iota$ from \eqref{iota} is given explicitly by
 $$\iota=\frac{1}{\sqrt 2\sqrt[4]{3}}\qmatrix{1&-1\\ \sqrt 3&\sqrt 3}\qquad \iota^{-1} = \frac{1}{\sqrt{2}\sqrt[4]{3}}\begin{pmatrix}
 \sqrt 3 & 1\\ -\sqrt 3 & 1
 \end{pmatrix}.
$$
\subsection{The Quantum space {of the torus phase space}}\label{qspace}
We recall the Weyl—Heisenberg operators, 
\begin{align}\label{WH}
(T^h_{\binom{y^*}{0}}\psi)(y)=\psi(y-y^*)&\qquad (T^h_{\binom{0}{\eta^*}}\psi)(y)=e^{\frac{i\eta^*}{h}y}\psi(y) \\
(T^h_{\binom{y^*}{\eta^*}}\psi)(y)&=\exp(\frac{i y^*\eta^*}{2h}) T^h_{\binom{y^*}{0}}T^h_{\binom{0}{\eta^*}},\nonumber
\end{align}

{and define a family of subspaces of $\mathcal{S}'(\RR)$ in terms of them,
\eq{
\hilb_{h,\ka}=\left\{\psi\in\ca{S}'(\RR): T_{\binom10}^h \psi= e^{2\pi i\ka_1}\psi,\ T_{\binom01}^h\psi=e^{2\pi i\ka_2}\psi\right\},\qquad \kappa\in[0,1]^2.
}{}
From \cite[Proposition 2.1]{BD96}, this space is non trivial if and only if $h=\frac{1}{2\pi N}$ for $N\in\NN^*$. Denoting the non-trivial spaces by $\mathcal{H}_{N,\ka},$ each of them is an $N-$dimensional Hilbert space {spanned by $\left\{\psi_{j,N}:=\sum_{n\in\ZZ}e^{-2\pi i\ka_1 n}\dl_{\frac {j+\ka_2}{N}+n},j=0,\dots,N-1\right\}$}.
{ Let us identify $\mathcal{H}_{N,\ka}$ and $\CC^N$ through the {(non-canonical)} map
\eq{
\psi = \sum_{j=0}^{N-1}\sum_{n\in\ZZ} \psi_j e^{-2\pi \ka_1 n} \dl_{\frac{j+\ka_2}{N}+n} \mapsto (\psi_0,\dots, \psi_{N-1})
}{iden} }

and equip $\mathcal{H}_{N,\ka}$ with the {(non-standard)} scalar product
\eq{
\inp{\phi_N,\psi_N}_{\mathcal{H}_{N,\ka}}=\frac 1N\sum_{j=1}^N \phi_{N,j}\overline{\psi_{N,j}},\qquad\qquad \phi_N=(\phi_{N,j})_{j=0}^{N-1},\psi_N=(\psi_{N,j})_{j=0}^{N-1}.
}{his}}

\begin{rem}
	Note that our choice of basis and inner product differs from the one in \cite[Proposition 2.1, Proposition 2.3]{BD96}, in which the authors chose as basis $\{\frac{1}{\sqrt N}\psi_{j,N}\}_{j=0}^{N-1}$. This change of basis reflects on the inner product \eqref{his} by multiplication in a factor\footnote{which influences our choice of norm, cf. \cite[Proposition 2.3.(ii)]{BD96}} $N$. However, our results \autoref{main},\autoref{prop:blabla} do not depend on that choice (and an analogue of \autoref{lacor} holds as well when using the conventions of \cite{BD96}).
\end{rem}

{One can decompose} $L^2(\RR)$ as a direct integral of these spaces (cf. \cite[Proposition 2.3.iii]{BD96}),

 $$L^2(\RR)=\int_{[0,1]^2}^\oplus \hilb_{N,\ka} d\ka.$$ 
 {For every $\ZZ^2-$periodic symbol $a\in C^\8(\RR^2)$ {(or equivalently $a\in C^\8(\TT^2)$)} $\Op_h(a):\Ltr\to\Ltr$ is bounded, acts on $\mathscr{S}'(\RR)$ and preserves each $\mathcal{H}_{N,\ka}$}, \eq{\Op_h(a)=\int^\oplus_{[0,1]^2} \Op_{N,\ka}(a) d\ka,}{dint}
{with $\Op_{N,\ka}(a):\mathcal{H}_{N,\ka}\to \mathcal{H}_{N,\ka}$ being the restriction of $\Op_h(a)$ to $\mathcal{H}_{N,\ka}$.} We recall the following lemma

\begin{lem}[Algebraic properties of $\Op_{N,\ka}(a)$, {\cite[Theorem XIII.83]{RS78}}] \label{lemnk1}
Let $a\in C^\8(\TT^2)$ real-valued, then
\begin{enumerate}[ref =\autoref{lemnk1}.\roman*]
\item\label{lemnksa} {For every $\phi_N,\psi_N\in\mathcal{H}_{N,\ka},$ $\inp{(\Op_{N,\ka}(a))\phi_N,\psi_N}=\inp{\phi_N,\Op_{N,\ka}(a)\psi_N}$.}
\item\label{lemnksup} $\norm{\Op_h(a)}{\mathcal{B}(L^2(\RR))}{}=\sup_{\ka\in [0,1]^2}\norm{\Op_{N,\ka}(a)}{\mathcal{B}(\hilb_{N,\ka})}{}$. \label{lemitem:2}
\end{enumerate} 
\end{lem}

\begin{rem}\label{remprodop}
	The proof of \ref{lemnksup} implies that the statement holds for finite products of pseudo-differential operators: Suppose $a_1,\dots, a_n\in C^\8(\TT^2)$ then  $$\norm{\prod_{j=1}^n\Op_h(a_j)}{\bltr}{}=\sup_{\ka\in[0,1]^2}\norm{\prod_{j=1}^n\Op_{N,\ka}(a_j)}{\mathcal{B}(\hilb_{N,\ka})}{}.$$
\end{rem}
For every $\gamma\in\tilde\Gamma(2)$ the map $\mathcal{M}_h(\gamma)$ sends $\hilb_{N,\ka}$ to $\hilb_{N,{}^t\gamma\ka}$ (cf. \cite{D93}), hence we will only focus on $\hilb_{N,0}$ for which the reduced operator $\mathcal{M}_{N,0}:=\mathcal M_h\mid_{\mathcal{H}_{N,0}}$, such that $\mathcal{M}_{N,0}:\hilb_{N,0}\to\hilb_{N,0}$. This space has the basis $\{e_j\}_{j=0}^{N-1}$ with $e_j=\dl_{\frac{j}{N}+\ZZ}$ and is viewed as a Hilbert space by \eqref{his}.

A general formula for $\mathcal{M}_{N,0}$ was given in \cite{HB80}: For every $\gamma=\psmallmatrix{a_{11}&a_{12}\\ a_{21}&a_{22}}\in\widetilde\Gamma(2)$ such that $a_{12}\neq 0$
(this is the case for all hyperbolic matrices) the elements of
$\mathcal{M}_{N,0}$ are then given by
\begin{equation}\label{metk}
	(\mathcal{M}_{N,0})_{q_2 q_1} = \frac{1}{\sqrt{i\abs{a_{12}}N}}\sum_{k=0}^{\abs{a_{12}}-1}{\exp\left(\frac{\pi i}{a_{12}N}\left(a_{11}(q_1+Nk)^2-2q_2(q_1+Nk)+a_{22}q_2^2\right)\right)}.
\end{equation}

One considers the stable and the unstable foliations on $\TT^2$:  The global leaf at the point{ $(y_1,\eta_1)\in\TT^2$ is defined by $W_{u,(y_1,\eta_1)}^{\TT^2}=\bigslant{W_{u,(y_0,\eta_0)}^{\RR^2}}{\ZZ^2}$ and $W_{s,(y_1,\eta_1)}^{\TT^2}=\bigslant{W_{s,(y_0,\eta_0)}^{\RR^2}}{\ZZ^2}$ for $(y_0,\eta_0)\equiv (y_1,\eta_1)\mod 1$.}
\begin{rem} For every $A=\Op_h(a)$ with $ a\in S_{L,\rho+}(\TT^2)$  we denote for brevity 
	\[A_N=\Op_{N,0}(a)=\Op_N(a)\qquad \mathcal{M}_{N,0}=\met\qquad \hi=\hilb_{N,0}.\] With $\gamma$ defined above we denote the eigenfunctions of $\met(\gamma)$ by $\{\vphi_{j,N}\}$.
\end{rem} 
Restricting to a symbol $a\in S_{L,\rho+}(\TT^2)$ analogous properties {(proved for sake of completeness)}
\begin{lem}[Analytical properties of $\Op_{N}(a)$]\label{lemnk2}
Let $a\in S_{L,\rho+}(\TT^2)$, then
\begin{enumerate}[ref =\autoref{lemnk2}.\roman*]
\item\label{lembd} $\Op_{N}(a):\mathcal{H}_{N}\to\mathcal{H}_{N}$ is bounded uniformly in $N$.
\item\label{lemnkmoy} A Moyal product formula holds: for $d\in S_{L,\rho+}(\TT^2)$, $\Op_{N}(a)\Op_{N}(d)=\Op_{N}(c)$ with $c=ad+O_{S_{L,\rho+}(\TT^2)}(h^{1-\rho-})$.
\end{enumerate} 
\end{lem}

\begin{proof}
	\begin{enumerate}[leftmargin=*]
		\item[]
		\item Follows immediately from combining \ref{lemnksup} with the boundedness of $\Op_h(a)$ proved in \autoref{bdopa} (potentially after rescalling).
		\item From \ref{slr+:moyal}, a Moyal product formula holds when considering the Weyl quantization of the symbols,
		\(
		\Op_h(a)\Op_h(d) = \Op_h(c)
		\)
		with $c=ad + O(h^{1-\rho-})$. Since $L^2(\RR)$ and the quantization split we can write
		\begin{align*}
			\int^\oplus_{[0,1]^2} \Op_{N,\ka}(a)\Op_{N,\ka}(d)d\ka = \Op_h(a)\Op_h(d) = \Op_h(ad) + \Op_h(r) = \int^\oplus_{[0,1]^2} \Op_{N,\ka} (ad) + \Op_{N,\ka}(r)d\ka
		\end{align*}
		with $\norm{\Op_h(r)}{\bltr}{}=O(h^{1-\rho-})$. Then $$\Phi_{\phi}(\psi) = \inp{(\Op_{h}(ad)+\Op_{h}(r) - \Op_{h}(a)\Op_h(d))\phi,\psi}$$ is the zero functional in $(\Ltr)^*$ and its projection to any $\ka-$fiber is the zero functional in $(\mathcal{H}_{N,\ka})^*$. From \ref{lemnksup}, $\norm{\Op_{N,\ka}(r)}{\ca{B}(\hilb_{N,\ka})}{} \le \norm{\Op_{h}(r)}{\mathcal{B}(L^2(\RR))}{}=O(h^{1-\rho-})$.
	\end{enumerate}
\end{proof}

{Similar analytical properties hold when $a\in S_{L,\rho,0}(\TT^2)$. In fact,}
\begin{lem}[Moyal product for symbols in $S_{L,\rho,0}(\TT^2)$]\label{moy0}
	Let $a,b\in S_{L,\rho,0}(\TT^2)$ then $$\Op_{N}(a)\Op_{N}(b) = \Op_{N}(c)$$ where $c=ab+O_{S_{L,\rho,0}(\TT^2)}(h^{1-\rho})$.
\end{lem}
\begin{rem}\label{nsp}
{Both in \ref{lemnkmoy} and in \autoref{moy0} if the supports of $a$ and $b$ are disjoint, repeated integration by parts yields that $\Op_{N}(c)=O_{\mathcal{B}(\hi)}(N^{-\8})$}
\end{rem}

We formulate a version of G\r arding inequalities holding for $\Op_{N}(a)$.
\begin{lem}\label{lemga}
 Let $a\in S_{L^s,\rho,\rho'}(\RR^2)$. Suppose $a\in S_{L^s,\rho,\rho'}(\TT^2)$ and $a\ge 0$ then there is some $C>0$ and $h_0>0$ such that for $h\in (0,h_0)$ $$\inp{\Op_{N}(a)u,u}\ge -Ch^{1-\rho-\rho'}\norm{u}{\hilb_{N,0}}{2}\qquad u\in\hilb_{N,0}.$$
\end{lem} 
{
\begin{proof}
Let $a\in S_{L^s,\rho,\rho'}(\RR^2)$ such that $a\ge 0$. We deduce from G\r arding inequality for $S_{L^s,\rho,\rho'}(\RR^2)$ that the spectrum of $\Op_h(a)$ lies inside the band $[-Ch^{1-\rho-\rho'},\norm{\Op_h(a)}{}{}]$ and
\[
-Ch^{1-\rho-\rho'}\text{Id}_{\Ltr} \le \Op_h(a) \le \norm{\Op_h(a)}{}{}\text{Id}_{\Ltr}
\] Denoting $\mathscr{a}=a-\ha\norm{\Op_h(a)}{}{}$ the last equation reads
\[
-\ha(\norm{\Op_h(a)}{}{} + Ch^{1-\rho-\rho'})\text{Id}_{\Ltr} \le \Op_h(\mathscr a) \le\ha(\norm{\Op_h(a)}{}{} + Ch^{1-\rho-\rho'})\text{Id}_{\Ltr}
\]
hence the self-adjointness of $\Op_h(a)$ implies $\norm{\Op_h(\mathscr a)}{}{} \le \ha(\norm{\Op_h(a)}{}{} + Ch^{1-\rho-\rho'})$. Then from \ref{lemnksup} also
\[
\norm{\Op_{N}(\mathscr a)}{}{} \le \ha(\norm{\Op_h(a)}{}{} + Ch^{1-\rho-\rho'})
\]
and from the self-adjointness of $\Op_{N}(a)$
\[
-\ha(\norm{\Op_h(a)}{}{} + Ch^{1-\rho-\rho'})\text{Id}_{\mathcal{H}_{N,0}} \le \Op_N(\mathscr a) \le\ha(\norm{\Op_h(a)}{}{} + Ch^{1-\rho-\rho'})\text{Id}_{\mathcal{H}_{N,0}}
\]
implying
\[
\Op_N(a)\ge -Ch^{1-\rho-\rho'}\text{Id}_{\mathcal{H}_{N,0}}
\]
\end{proof}
}
We note that $\Op_{N}$ satisfies an exact Egorov relation with the dynamics (cf. \cite{DG03}). Namely {the restriction of \eqref{egorov_full} to $\hi$ means that} for every $\gamma\in \tilde \Gamma(2)$, 
\eq{
\met(\gamma)^*\Op_{N}(a)\met(\gamma)=\Op_{N}(a\circ \gamma),\qquad a\in C^\8(\TT^2).
}{egorov}
One can express $\met(\gamma)=(\met(\gamma)_{jk})_{j,k}\in \CC^{N\times N}$ explicitly by number theoretic sums. For instance for $\gamma_{DE}$, for every $N\in 2\NN+1$, up to some scalar phase factor (its expression as a quotient of Gauss sums is expanded in \cite{KR00} and \cite{DG03}),
\[
\met(\gamma_{DE})_{jk} = \frac{1}{\sqrt N}\exp(\frac{2\pi i}{N}(k^2-kj+j^2)).
\]
{We note indeed that the quadratic phase of $\mathcal{M}_h(\gamma)$  evaluated at points $y,\eta\in \frac{1}{\NN}\ZZ$ is exactly the phase of $\met(\gamma)$.}

\section{{The main estimate}}\label{quanmain}
Let $f_N,g_N\in \hi,l=(l_1,l_2)\in\ZZ^2$. Recall the discrete 
Fourier—Wigner (FW) transform,{ \eq{\mathscr{V}_N(f_N,g_N)(l_1,l_2)=\inp{T_{\frac lN}^h f_N, g_N}_{\hi}=\frac {e^{\frac{i l_1l_2}{N}}}{N} \sum_{k=0}^{N-1} e^{\frac{2\pi i l_2k}{N}}f_{N,k-l_1}\overline{g_{N,k},}}{WT}
with $(f_{N,k-l_1})_{k=0}^{N-1}:=(f_{N-l_1},\dots,f_{N-1},f_0,\dots, f_{N-l_1-1})$ and $T_{\frac lN}^h:\hi\to\hi$ being the restriction of the Weyl—Heisenberg operator in \eqref{WH} to $\hi$. { In order to pass from $\Op_h(a)$ to $\Op_{N}(a)$ one recalls that $h=\frac{1}{2\pi N}$ and restrics $(y^*,\eta^*)$ to lie in $(2\pi \ZZ)^2$.}
	 We will study the localization of eigenfunctions of $\met(\gamma)$ using $\mathscr{V}_N$, or more precisely { its Fourier transform, which is} the discrete Wigner distribution.

Consider the diagonal matrix coefficients $\inp{\Op_N(a)\vphi_{j,N},\vphi_{j,N}}$. 
\begin{defn}
{The distribution $\mathscr{W}_{j,N}$ defined by $\inp{a,\mathscr{W}_{j,N}}_{\mathcal{D}'(\TT^2)} := \inp{\Op_N(a)\vphi_{j,N},\vphi_{j,N}}$ is called the diagonal Wigner distribution.} It can be expressed explicitly in terms of $T_{\frac lN}^h$ (cf. \cite{BD96}). Denoting the Fourier coefficients of $a$ by $\check a:=\{\hat a(l)\}_{l\in\ZZ^2}$,
\begin{align}\label{Wad}
 \inp{\Op_N(a)\vphi_{j,N},\vphi_{j,N}} &= \inp{\sum_{l\in\ZZ^2} \hat a(l)T_{\frac lN}^h\vphi_{j,N},\vphi_{j,N}}_{\hi} \\
&= \sum_{(l_1,l_2)\in\ZZ^2}\int_{\TT^2} a(x,\xi;h)e^{-2\pi i(x l_2-\xi l_1)}dxd\xi \cdot \mathscr{V}_N(\vphi_{j,N},\vphi_{j,N})(l_1,l_2)\nonumber \\
&= \inp{\check a,\check{\mathscr{V}_N}(\vphi_{j,N},\vphi_{j,N})}_{\CC^{\ZZ^2}} \nonumber 
\end{align}
{with $\check{\mathscr{V}_N}(\vphi_{j,N},\vphi_{j,N}) := \{\mathscr{V}_N(\vphi_{j,N},\vphi_{j,N})(l)\}_{l\in\ZZ^2}$ and $\check a(l)$ being the $l-$th symplectic Fourier coefficient, i.e,.
	\[
	\check a(l)=\int_{\TT^2} a(x,\xi;h)e^{-2\pi i(x l_2-\xi l_1)}dxd\xi.
	\]
Since $\mathscr{V}_N$ has moderate growth and the Fourier coefficients of $a$ are exponentially decaying this expression is well defined thus applying Plancherel's identity we arrive to a distribution $\mathscr{W}_{j,N}\in\mathscr{D}'(\TT^2)$ having Fourier coefficients $$\{\widehat{\mathscr{W}_{j,N}}(l)\}_{l\in\ZZ^2}:=\check{\mathscr{V}_N}(\vphi_{j,N},\vphi_{j,N})$$.}
\end{defn}  
{In fact \cite{HB80} and \cite{BD96} show that $\mathscr{W}_{j,N}$ is a linear combination of Dirac peaks on a lattice.} One can extract {from $\{\vphi_{j,N}\}$ a sub-sequence $\{\vphi_{j_{k},N_k}\}$ such that $\mathscr{W}_{j_k,N_k}\to \mu_{\text{sc}}$} as $N_k\to\8$. That is done in \cite{BD96} by considering a positive quantization {replacing $\mathscr{W}_{j,N}$ by probability measures}. {Every such limit} $\mu_{\text{sc}}$ is called {a semi-classical measure} and we say $\{\mathscr{W}_{j_{k},N_k}\}$ converges semi-classically to the measure $\mu_{\text{sc}}$. 

Rather than studying $\mathscr{W}_{j,N}$ directly we give bounds on an asymptotic expression, $$\inp{\Op_N(a)\vphi_{j,N},\Op_N(a)\vphi_{j,N}}=\inp{\Op_N(\abs{a}^2)\vphi_{j,N},\vphi_{j,N}}+O(N^{-1}),\qquad a\in C^\8(\TT^2).$$ We obtain estimates from below on the {left} hand side. In turn they yield a lower bound on {$\inp{a^2,\mathscr{W}_{j,N}}$} and in fact on the limit of its convergent sub-sequences. We adapt the methods introduced in \cite{DJ17}, acquiring a quantitative lower bound on $\norm{\Op_N(a)u}{}{}$ for every $u\in\hi$.
\begin{thm}[The main theorem]\label{prop:blabla}
Let $0\not\equiv a\in C^\8(\TT^2)$ and $\gamma\in\tilde\Gamma(2)$. Then there exist constants {$C_1(a),C_2(a),N(a)$ depending only on the choice of $a$ such that for every $N\ge N(a)$} and $u\in \hi$,
\[
\norm{u}{\hi}{} \le  
C_1(a) \norm{\Op_{N}(a)u}{\hi}{} + C_2(a)\log N\min_{z\in\CC:|z|=1}\norm{(\met(\gamma)-z)u}{\hi}{},
\]
\end{thm}

From this estimate we deduce \autoref{main}
\begin{proof}[\hypertarget{pff}{\proofname { of \autoref{main}}}]
Fix an open set $\emptyset\neq\Omg\subseteq\TT^2$. Let $\{\vphi_{N_j}\}_{j\to\8}$ be a subsequence of eigenfunctions of $\mathcal{M}_{N_j}(\gamma)$ converging semi-classically to $\mu_{\text{sc}}$. Fix a non-vanishing symbol $a\in C^\8(\TT^2)$ {supported inside  $\Omg$}. Applying \autoref{prop:blabla}, there exists a constant $C_1(a)\in\RR$ such that $\norm{\Op_N(a) \vphi_{j,N}}{}{} \ge C_1^{-1}(a)$ for $N$ large enough. From the semiclassical convergence $$\norm{\Op_N(a)\vphi_{j,N}}{}{2}\xrightarrow{N\to \8} \int_{\Omg}|a|^2(x,\xi)d\mu_{\text{sc}},$$ thus $\int_{\Omg}|a|^2d\mu_{\text{sc}}\ge C_1(a)^{-2}>0$ and { $\mu_{\text{sc}}(\Omg)\ge C_\Omg= C_1(a)^{-2}(\max_{\Omg}|a|)^{-1} >0$.}
\end{proof}

{Another result of the theorem is \autoref{lacor}:
\begin{proof}[\hypertarget{pflacor}{\proofname} of \autoref{lacor}]
	Fix $0\le \al_1 < \al_2 \le 1$. {Let us take a smoothed characteristic function $a=a(x)\in C^\8(\TT^2)$ supported {inside} of $[\al_1,\al_2]\times \TT$}. Fix an eigenstate  of $\met(\gamma)$, $\ \vphi_N=\sum_k\vphi_{N,k} \dl_{\frac kN+\ZZ}$. We note that in this scenario we can write explicitly $[\Op_N(a)\vphi_N] =\sum_{k\in \llbracket\al_1N,\al_2N\rrbracket} a(\frac kN)\vphi_{N,k}\dl_{\frac kN+\ZZ}$. From \autoref{prop:blabla} for these choices of a symbol and an eigenstate there exists a $N_0\ge {N(a)}\in\NN^*$ such that for every $N>N_0$
	\begin{align*}
	1=\norm{\vphi_N}{\hi}{2} &\le C_1^2(a)\norm{\Op_N(a)\vphi_N}{}{2}
	\end{align*}
	meaning 
	\[
	1  \le {\frac {C_1^2(a)}{N}}\sum_{k\in \llbracket\al_1N,\al_2N\rrbracket} \abs{a(\frac kN)\vphi_{N,k}}^2 \le {\frac {C_1^2(a)}{N}}\sum_{k\in \llbracket\al_1N,\al_2N\rrbracket} \abs{\vphi_{N,k}}^2.	
	\]
	\end{proof}
}
\subsection{Geometric construction and propagated operators}\label{gcon}
Fix $a\in C^\8(\TT^2)$ having a non-empty support and consider two proper non-empty open subsets $\mathcal{K}_1,\mathcal{K}_2\sub\TT^2$ with $\mathcal{K}_2\sub\supp(a)$. Let $a_1,a_2\in C^\8(\mathbb{T}^2)$ be a couple of symbols satisfying {(see \autoref{fig:pt})}
\begin{align}\label{geoconstr}
 0\le a_\bt \le 1&, \quad a_1+a_2=1,\quad \supp(a_1)\sub(\supp(a))^\circ \quad a_1\mid_{\mathcal{K}_2}\equiv a_2\mid_{\mathcal{K}_1} \equiv 1\nonumber
\end{align}
 with the corresponding operators $\Op_N(a_1),\Op_N(a_2)$ satisfying  \eq{\tag{$\star$} A_{1,N}+A_{2,N} = \text{Id}_{\hi}\qquad A_{\bt,N}=\Op_N(a_\bt).}{streq}\stepcounter{equation}
\begin{figure}[tp]
\centering
\includegraphics[scale=0.4]{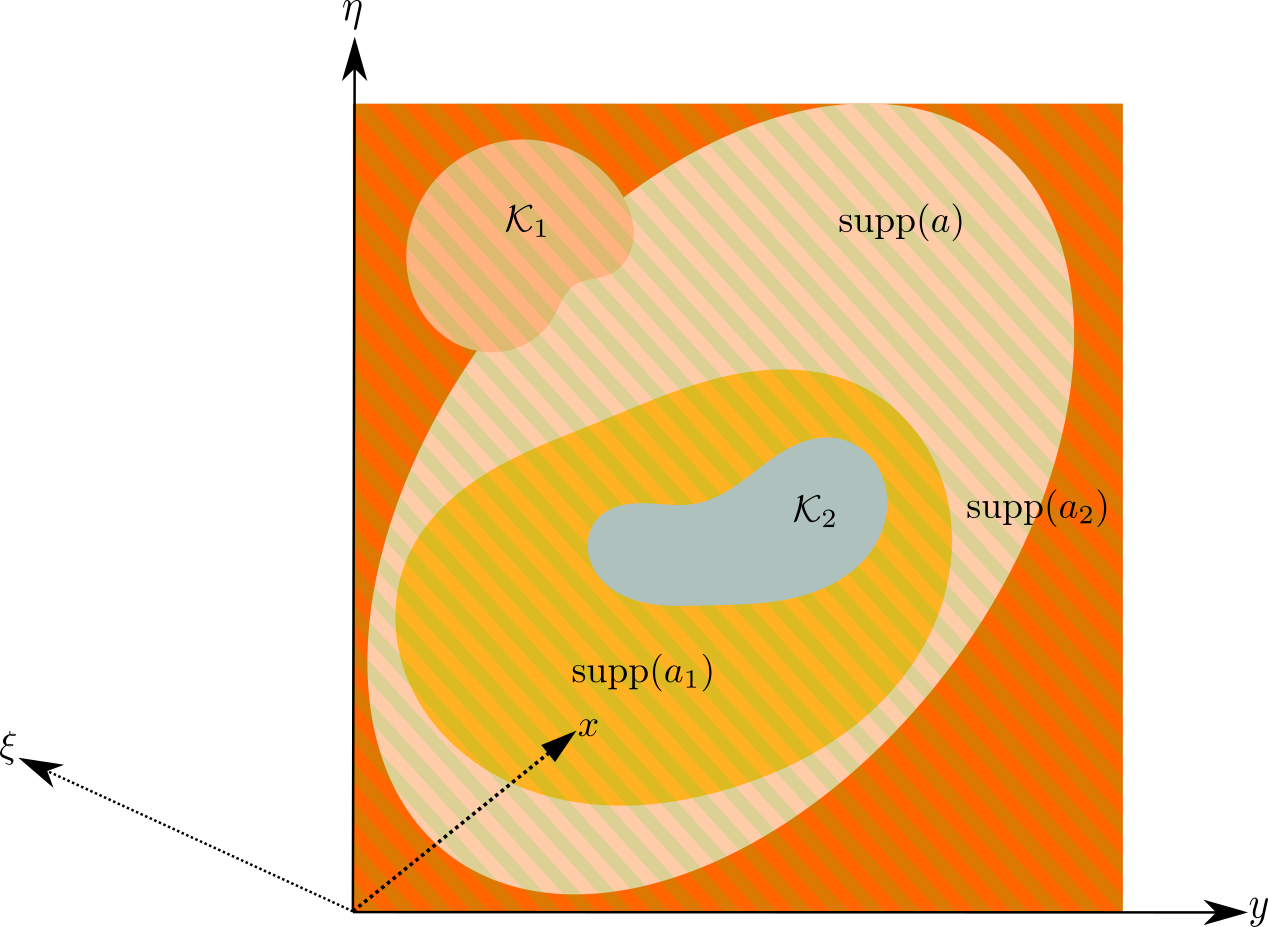}
\caption{{ Given a symbol $a\in C^\8(\TT^2)$ we construct in \autoref{gcon} a partition of unity on $\TT^2\cong \left[0,1\right)^2$. $\RR^2$ is equipped with $(y,\eta)$ coordinates. The $(x,\xi)$ coordinates are obtained after applying the coordinate map $\iota$ to the coordinates $(y,\eta)$. The new coordinates describe the decomposition of $\RR^2$ into stable and unstable directions.} The yellow domain represents $\supp(a_1)$ and the domain in green lines represents $\supp(a_2)$.} 
\label{fig:pt}
\end{figure}
We can lift our construction to $\RR^2$ by periodization: Pick $L^{u}_{\RR^2}$ to be the Lagrangian foliation whose projection is $L^u$. {We can lift $a_\bt$ to $\RR^2$}. Quantizing these symbols yields the corresponding pseudo-differential operators $ A_\bt=\Op_h(a_\bt):\Ltr\to\Ltr$ satisfying \stepcounter{equation}
\eq{
{a_1+a_2=1,\quad \Op_h(a_1)+\Op_h(a_2)=\text{Id}_{\Ltr}, \quad a_1\mid_{\mathcal{K}_2+\ZZ^2}\equiv a_2\mid_{\mathcal{K}_1+\ZZ^2} \equiv1.}
\tag{$\br\star$}}{starbar}
{ By an abuse of notation henceforth we identify {$\ZZ^2$-periodic} symbols $a:\RR^2\to\RR$ with their restrictions to $\TT^2$}. Henceforth we fix some $\rho\in (\ha,1)$. Given a pseudo-differential operator $A=\Op_h(a):L^2(\RR)\to L^2(\RR)$ we use the unitarity of $\met(\gamma)$ in order to propagate with respect to the hyperbolic dynamics up to time $n=0,\dots,T'$ with 
\eq{ T = \left\lfloor \frac{\rho}{4\log|\lm_u|}\log\frac 1h\right\rfloor,\quad T'=4T.\tag{T}}{T}\stepcounter{equation} Let us denote $ A(n)=\ca{M}_h^{-n}(\gamma)A\ca{M}_h^n(\gamma)$. Recalling the exact Egorov relation \eqref{egorov_full}, $ A(n)=\Op_h(a\circ\gamma^n)$. For the restriction $A_N:\hi\to\hi$ from applying the exact Egorov relation in \eqref{egorov}, $A_N(n)=\met(\gamma)^{-n}A_N\met(\gamma)^n=\Op_N(a\circ\gamma^n)$. We note that by \autoref{lemLLxi} $a\circ\gamma^n\in S_{L^s,\rho+}(\TT^2)$ uniformly in $0\le n < T'$.

The partition of unity \eqref{starbar} on $\Ltr$ {together with the evolution through $M_h(\gamma)$ allows to construct refined partitions of unity, in terms of "word operators"}: {to} every word $w\in\mathcal{W}(n)=\{1,2\}^n=\{w=w_{0}\dots w_{n-1} : w_{j}\in\{1,2\}\}$ {(where $0 \le n \le T'$)} corresponds an operator on $\Ltr$ \begin{align*}{A_w}&= A_{w_{n-1}}(n-1)A_{w_{n-2}}(n-2)\cdots A_{w_1}(1)A_{w_0}\\&=\Op_h(a_{w_{n-1}}\circ \gamma^{n-1})\cdots  \Op_h(a_{w_1}\circ\gamma)\Op_h(a_{w_0})\end{align*}
{ as well as an operator on $\hi$, $A_{w,N}=(A_w)_N$. Suppose $B=\Op_h(b)$ with $b\in S_{L^s,\frac\rho4,0}(\RR^2)$ we denote the "tilted" operator $\widetilde{B}=\Op_h(\tilde b)=\ca{M}_h^{-1}(\iota) B\ca{M}_h(\iota)$ having from \autoref{lemLLxi} the "tilted" symbol $\tilde b\in S_{L^\xi,\frac\rho4,0}(\RR^2)$.}

\begin{lem}\label{lemaw}
	\begin{enumerate}[leftmargin=*, ref=\autoref{lemaw}.\roman*]
		\item[]
		\item Let $0\le n\le 4T$ then $a_w\in S_{L^s,\rho+}(\RR^2)$ and $A_w=\Op_h(a_w)+O_{\cal{B}(\Ltr)}(h^{1-\rho-})$ where $a_w=\prod_{j=0}^{n-1}a_{w_j}\circ\gamma^j$ uniformly in $n$ and in $w\in\mathcal{W}(n)$. 
		\item\label{lemawiii} {If $0\le n\le T$ then $a_w\in S_{L^s,\frac{\rho}{4}+}(\RR^2)$ and $A_w=\Op_h(a_w)+O_{\cal{B}(\Ltr)}(h^{1-\frac{\rho}{4}-})$ where $a_w=\prod_{j=0}^{n-1}a_{w_j}\circ\gamma^j$ uniformly in $n$ and in $w\in\mathcal{W}(n)$. }
	\end{enumerate}
\end{lem}

\begin{proof}
We prove only the first statement in the lemma as (iii) follows by replacing $\rho$ by $\frac\rho4$ and the others by exchanging the stable and unstable directions. Consider the tilted operator $\widetilde{A_w}$ which from applying \eqref{egorov_full} can be expressed as
\begin{align*}
\widetilde{A_w}
= \Op_h(\stackon[-8pt]{a_{w_{n-1}}\circ \gamma^{n-1}}{\vstretch{1.5}{\hstretch{2.4}{\widetilde{\phantom{\;\;\;\;\;\;\;\;}}}}}) \cdots  \Op_h(\widetilde{ a_{w_1}\circ\gamma})\Op_h(\widetilde{a_{w_0}}).
\end{align*}
{Let us now} show that $\widetilde{A_w}$ can be expressed as $\widetilde{A_w}=\Op_h(\widetilde{a_w})+O(h^{1-\rho-})$
where \begin{align}\label{pra}
	\widetilde{a_w}=\prod_{j=0}^{n-1}b_j :=\prod_{j=0}^{n-1}a_{w_j}\circ\gamma^j\circ\iota\in S_{L^\xi,\rho+}(\RR^2).
	\end{align}
We note that {from applying Leibniz's chain rule} $\bnd_y^\al\bnd_\eta^\be \widetilde{a_w}$ is a sum of $n^{\al+\be}=O(h^{-})$ summands, each of the form \eq{d(y,\eta)=(\bnd_y^{\al_{0}}\bnd_\eta^{\be_{0}}b_0)\cdots (\bnd_y^{\al_{n-1}}\bnd_\eta^{\be_{n-1}}b_{n-1})(y,\eta)}{theprod} with $\sum_j\al_{j}=\al,\sum_k \be_{k}=\be$. By \autoref{lemLLxi} $b_{j}\in S_{L^\xi,\rho,0}(\RR^2)$ hence each factor in the product \eqref{theprod} is bounded by $C_{\al_j,\be_j}h^{-\rho\al_j}$ {uniformly in $w$ and in $j$}. However, since the number of derivatives, with respect either to $y$ or to $\eta$ does not exceed $\al+\be$ and since $|b_j|\le 1$,  
\eq{
|d(y,\eta)| \le \prod_{j:(\al_j,\be_j)\neq 0}C_{\al_j,\be_j} h^{-\rho \sum \al_j}\le \max_{j} C_{\al_j,\be_j}^{\al+\be} \cdot h^{-\rho \al}.
}{dbnd}
{We remind $\al_j\in\llbracket0,\al\rrbracket,\be_j\in\llbracket0,\be\rrbracket$ and that $C_{\al_j,\be_j}$ depend only on the number of derivatives taken. As follows, for every $\al,\be,\ $ $\max_{j}C_{\al_j,\be_j}^{\al+\be}$ is uniformly bounded by some constant $C$. Since there are $O(h^{-\eps})$ terms, each of them bounded by some $Ch^{-\al\rho}$ we deduce that $\abs{\bnd_y^\al\bnd_\eta^\be \widetilde{a_w}} \le Ch^{-\al\rho-}$ and $\widetilde{a_w} \in S_{L^\xi,\rho+}(\RR^2)$}

We prove that uniformly with respect to the word $w\in\mathcal{W}(n)$ where $1\le n\le 4T,$
\eq{
\widetilde{A_w} - \Op_h(\widetilde{a_w}) = O_{\bltr}(h^{1-\rho-}).
}{remnd}
{It can be verified by induction} on $2\le n\le 4T$ that
\eq{
\widetilde{A_w}-\Op_h(\widetilde{a_w}) = \sum_{k=1}^{n-1} B_k \Op_h({b_{w_{k+1}}})\cdots \Op_h({b_{w_{n-1}}}) ,
}{di4}

where $$B_k=
 \Op_h(b_{k})\Op_h(\prod_{j=0}^{k-1} b_{j})-\Op_h(\prod_{j=0}^{k} b_j).$$ 

We remind that $\prod_{k=0}^{j-2} b_{k}\in S_{L^\xi,\rho+}$ uniformly with respect to $j\in\llbracket0,n-1\rrbracket$. Let us show a bound on the norm of this operator: 
{
We first note that for every $1\le j\le 4T-1$, from the Moyal product in $S_{L^\xi,\rho,+}(\RR^2)$ (from \ref{slr+:moyal}) we deduce that $\norm{B_j}{\bltr}{} = O(h^{1-\rho-})$ uniformly in $j$. Then following \autoref{bdopa} we deduce the estimate 
\begin{align*}
\norm{\Op_h({b_{{n-1}}})\cdots \Op_h({b_{{j+1}}}) B_j}{\bltr}{} &\le C(1+ C'' h^{1-\rho})^{C'\log\frac 1h}h^{1-\rho-}
\end{align*}
for constants $C,C''$ independant of $j$ and $n$.
{As $h\to 0$
\[
\log (1+ C_k h^{1-\rho})^{C'\log\frac 1h}  = -C'\log h\cdot\log(1+C_k h^{1-\rho}) = O(h^{1-\rho}\log h) \to 0
\]
}
thus $\frac{\norm{\Op_h({b_{w_{n-1}}})\cdots \Op_h({b_{w_{k+1}}}) B_k}{}{}}{h^{1-\rho-\eps}}$ remains uniformly bounded as $h\to 0$. As the sum is over $n$ operators of this form, we deduce from the triangle inequality that 
\[
\norm{\widetilde{A_w}-\Op_h(\widetilde{a_w})}{}{} \le C''h^{1-\rho-}.
\]
{With $C''$ being uniform with respect to $w$ and $1\le n\le 4T$.}
}
\end{proof}

{Let us now consider operators which correspond to shorter "words", of length $T$.} For every function $c:\mathcal{W}(T)\to\RR$ (with $T$ defined in \eqref{T}) such that $\sup|c(w)|\le 1$ corresponds the weighted operator $A_{c}=\sum_{w\in \mathcal{W}(T)}c(w)A_w$. It is a pseudo-differential operator with a principal symbol $a_c=\sum_{w\in \mathcal{W}(T)}c(w)a_w$ or equivalently, 

\begin{lem}\label{lemac}
Suppose $\sup_{w\in \mathcal{W}(T)}|c(w)|\le 1$ then $a_c\in {S_{L^s,\frac{1}{2},\frac{1}{4}}}(\RR^2)$ and $$A_{c}=\Op_h(a_c)+O(h^{1-\frac\rho4(1+\frac{\log 2}{\log |\lm_u|})-\eps}).$$
\end{lem}

\begin{proof}
We first observe that $a_c$ is bounded uniformly
\[
|a_c| = \abs{\sum_{w\in \mathcal{W}(T)}c(w)a_w} \le \sum_{w\in \mathcal{W}(T)}|c(w)|a_w \le \sum_{w\in\mathcal{W}(T)} a_w = 1.
\]
{We estimate $\abs{\bnd_x^\al \bnd_\xi^\be \widetilde{a_c}}$ for every $\al,\be\in\NN$. First, by the second part of \autoref{lemaw}}, $\widetilde{a_w}\in S_{L^\xi,\frac{\rho}{4}+}(\RR^2)$ and recalling \autoref{lemsym}, if $(\al,\be)\neq(0,0)$
\begin{align*}
\abs{\bnd_x^\al \bnd_\xi^\be  \widetilde{a_c}} &\le \sum_{w\in \mathcal{W}(T)}\abs{\bnd_x^\al \bnd_\xi^\be  \widetilde{a_w}}\\ &\le \sum_{w\in\mathcal{W}(T)} C_{\al,\be}h^{-\frac{\rho}{4}\al-\eps} \\&\le 2^{T}C_{\al,\be} h^{-\frac\rho 4\al-\eps}
\end{align*}
 Recalling \autoref{lemaw}, $C_{\al,\be}$ is uniform. We recall \autoref{mineig} implying $|\lm_u| \ge 2+\sqrt 3$, $2^T\approx h^{-\frac{\rho\log 2}{4\log|\lm_u|}} \le h^{-\frac{\rho\log 2}{4\log|2+\sqrt 3|}} \le h^{-\frac{\rho}{4}+\rho\eps_0}$ for some $\eps_0>0$ and {hence taking $\eps < \frac{1-\rho}{4}$
\begin{align*}
\abs{\bnd_x^\al \bnd_\xi^\be  \widetilde{a_c}} \le C_{\al,\be}h^{-\frac{\rho}{4}-\frac{\rho}{4}\al+\rho\eps_0-\eps}\le h^{-\rho\frac\al 2- \rho\frac\be 4}
\end{align*}}
If $\al=\be=0, \ |\bnd_x^\al\bnd_\xi^\be\widetilde{a_c}| = \abs{\widetilde{a_c}}\le 1$
We deduce that $\widetilde{a_c}\in {S_{L^\xi,\frac{1}{2},\frac{1}{4}}}(\RR^2)$ and therefore $a_c\in {S_{L^s,\frac{1}{2},\frac{1}{4}}}(\RR^2)$. The remainder follows from applying \ref{lemawiii},
\[
A_c=\sum_{w\in\cal{W}(T)} c(w) A_w = \sum_{w\in\cal{W}(T)} c(w)(\Op_h(a_w)+O(h^{1-\frac\rho 4-\eps})).
\] 
The sum is over $2^T$ terms thus the reminder is of order $2^Th^{1-\frac\rho4-\eps} \le h^{1-\frac\rho 4(1+\frac{\log 2}{\log\abs{\lm_u}})-\eps}$
\end{proof}
{Henceforth for every subset $S\sub \mathcal{W}(T)$ we denote $A_{S,N}:=(A_{\ind_S})_N$.}
Fix { a small } $\delta\in(0,1)$ { that will be specified later in \eqref{ladel}}, and define the control function 
\eq{
F:\mathcal{W}(T)\to\RR,\quad F(w)=\frac{\#\{j\in\{0,\dots,T-1\}: w_j=1\}}{T}.\tag{$\text{F}$}}{F} Roughly speaking, {in the notations of \eqref{streq},} $F$ quantifies the fraction of time a word $w\in\mathcal{W}(T)$ "spends inside" {$\supp(a_1)$}. Therefore the set 
\[
\mathcal{Z}=\{w\in\mathcal{W}(T):F(w)\ge \delta \}
\]
{contains the "controlled short logarithmic" words, i.e., those who have corresponding symbols spending a fraction $\dl$ of time supported inside $\supp(a_1)$}. This definition opens the door for obtaining a {finer} pseudo-differential partition of unity. We would like to apply the fractal uncertainty principle (see \autoref{finalfups} below) hence we {consider a partition arising from} $\mathcal{W}(8T)$ (with $T$ from \eqref{T}), {words twice longer than Ehrenfest time $\frac{\log\frac 1h}{\log\abs{\lm_u}}\approx T'$}. The set $\mathcal{W}(8T)$ can be divided according to the control function $F$ to $\mathcal{W}(8T)=\mathcal{X}\sqcup\igrk$ with
\begin{align*}
    \mathcal{X} &= \{w^{(1)}w^{(2)}\cdots w^{(8)}\in\mathcal{W}(8T): w^{(j)}\in \mathcal{W}(T)\setminus \mathcal{Z},\ \forall j \} \\
    \mathcal{Y} &= \{w^{(1)}w^{(2)}\cdots w^{(8)}\in\mathcal{W}(8T): \forall k\ w^{(k)}\in\mathcal{W}(T)\text{ and }\exists j,\ w^{(j)}\in \mathcal{Z} \}.
\end{align*}
$\mathcal{Y}$ is the set of {"long logarithmic controlled words": words such that a single subword $w ^{(j)}$ of them is controlled.} $\mathcal{X}$ is the set of uncontrolled words: the time they spend inside this fixed subset of $\supp(a_1)$ can't be bounded from below. 
\subsection{Outline of the proof}

The proof of \autoref{prop:blabla} is based on two key estimates: First we estimate {the mass $\norm{A_{\igrk,N} u}{}{}$, i.e., coming from "controlled" words},

\begin{prop}\label{propigrk}
There are constants $C,C_1,C_2\in\RR$ satisfying for every {$N\ge 1$} and $u\in\hi$,
\[ 
\norm{A_{\igrk,N}u}{\hi}{} \le c_\dl(\gamma)(C_1\norm{\Op_{N}(a)u}{}{}+\log N\min_{z\in\CC:|z|=1}\norm{(\mathcal{M}_N(\gamma)-z)u}{\hilb_{N}}{}+\frac{C_2}{N^{\frac18}}\norm{u}{\hi}{})
\]
{with $c_\dl(\gamma)=\frac{C}{\dl}$.}
\end{prop}
The proof of this proposition is done in \autoref{seci} below and is relying on the almost-monotonicity property of the calculus {in \autoref{lemono}.} Then we turn to estimate the norm of operators which corresponds to uncontrolled words $w\in\mathcal{X}$. We estimate the norm of each $A_{w,N}u$ (for every $w\in\mathcal{X}$) separately by {an argument involving} the fractal uncertainty principle from \cite{BD16},
\begin{prop}\label{propfup}
There exists $\beta'>0$, a constant $C>0$ such that for every $N\in\NN^*$ $$\sup_{w\in\mathcal{W}(8T)}\norm{A_{w,N}}{\mathcal{B}(\hi)}{} \le \frac{C}{N^{{\be'}}},$$ with $T$ from \eqref{T}.
\end{prop}

Then, we invoke a combinatorial argument from \cite{DJ17},
\begin{lem}[Lemma 3.3 in \cite{DJ17}]\label{comb}
There is some $C$ (that might depend on $\delta$) such that for small enough $\dl$  {and for every $h,$} $\#\mathcal{X}\le Ch^{-4\sqrt\dl}$.
\end{lem}

Note that $N_0$ appearing in the proof of \cite[Lemma 3.3]{DJ17} is larger than $T$ specified in \eqref{T} above, yet one can deduce \autoref{comb} from \cite[Lemma 3.3]{DJ17}.

\begin{proof}[End of the proof of \autoref{prop:blabla}] Decomposing $\text{Id}_{\hi}=A_{\mathcal{X},N}+A_{\igrk,N}$ and applying the triangle inequality
\[
\norm{u}{\hi}{} \le \norm{A_{\cal X,N}u}{}{} + \norm{A_{\igrk,N} u}{}{}.
\]

From \autoref{comb}, \ref{lemnksup} and \autoref{propfup}, $\norm{A_{\cal X,N} u}{}{} \le C N^{-\be'+4\sqrt\dl}\norm{u}{}{}$. Combined with \autoref{propigrk}, it yields
\[
\norm{u}{}{}\le CN^{-\be'+4\sqrt\dl}\norm{u}{}{}+ \frac{C_1}{\dl}\norm{\Op_N(a)u}{}{}+\frac{C_2}{\dl}\log N\min_{z\in\CC:|z|=1}\norm{(\met(\gamma)-z)u}{}{} + \frac{C_3}{N^{\frac 18}}\norm{u}{}{}.
\] 
Choosing $\dl$ such that \eq{\be'-4\sqrt\dl=\frac{\be'} 4}{ladel} and picking $N$ big enough {eliminates} the first term which in turn implies the result. 
\end{proof}

\subsection{An estimate on the controlled region $\mathcal{Y}$}\label{seci}
We first bound from above $\norm{\Op_N(a_1)}{}{}$.
\begin{lem}\label{esta1} There exists some $C_1'\in\RR$ such that for every $u\in\hi$
\[
\norm{\Op_N(a_1)u}{\hi}{} \le C_1'\norm{\Op_N(a)u}{\hi}{} {+O(N^{-\8})\norm{u}{\hi}{}}.
\]
\end{lem}

\begin{proof}

{{First we recall that $\supp(a_1)\sub(\supp(a))^\circ$ and thus the function $\frac{a_1}{a}$ is well defined and one can construct a parametrix recursively.} We construct a function $q\in C^\8(\TT^2)$ satisfying $q\#_h a=a_1+O_{C^\8(\TT^2)}(N^{-\8})$ for every $(x,\xi)\in(\supp(a))^\circ$. Expressing by a formal ansatz $q(x,\xi)=\sum_{j=0}^\8 q_jh^j$ and using the expansion of the Moyal product in the left hand side we obtain the values of sequence of functions $\{q_j\}_{j=0}^{\8}$. For every $j\in\NN$ one can write explicitly $q_j=\frac{p(a,a_1,q_0,\dots,q_{j-1})}{C_j a}$ for a constant $C_j\in\RR$ and $p\in\RR[x_1,\dots,x_{j+1}]$. {Inductively $a_1,q_0,\dots,q_{j-1}$ are all continuous functions supported inside $(\supp(a))^\circ$ and thus $q_j$ is a well defined function on $\TT^2$}. By Borel's lemma one can indeed construct a symbol $q\sim \sum q_jh^j$ and we obtain $$\Op_N(a_1)=\Op_N(q)\Op_N(a)+O(N^{-\8}). $$ We deduce 
\begin{align*}
\norm{\Op_N(a_1)u}{}{} &\le \norm{\Op_N(q)\Op_N(a)u}{}{}+O(N^{-\8})\norm{u}{}{}\\&\le C_1\norm{\Op_N(a)u}{}{}+O(N^{-\8})\norm{u}{}{},
\end{align*}
obtaining the lemma.}
\end{proof}

Next, controlling the growth of $A_Nu$ implies control on the growth of $A_N(m) u$ for $u\in\hi$,
\begin{lem}\label{estA}
Let {$A_N:\hi\to \hi$} be a pseudo-differential bounded operator, then for all $m\in\ZZ$,
\[
\norm{A_N(m) u}{\hi}{} \le \norm{A_N u}{\hi}{} +|m|\norm{A_N}{\cal{B}(\hi)}{}\min_{z\in\CC:|z|=1}\norm{(\met(\gamma)-z)u}{\hi}{},
\] where $u\in \hi$.
\end{lem}

This lemma is the discrete analogue of propagation of singularities for long time (cf. \cite[Theorem E.47]{DZ19}).

\begin{proof}
	{
Without loss of generality assume that $m\in\NN$. Let $z\in\CC$ having $|z|=1$. From the triangle inequality,
\begin{align*}
\norm{A_N(m)u}{\hi}{} = \norm{A_N\met^m(\gamma)u}{\hi}{} &\le \norm{A_Nu}{\hi}{} + \norm{(A_N\met^m(\gamma)-z^mA_N)u}{}{} \\
&\le \norm{A_Nu}{\hi}{} + \norm{A_N}{\mathcal{B}(\hi)}{}\norm{(\met^m(\gamma)-z^m\text{Id}_{\hi}) u}{}{}.
\end{align*}
Since $\met^m(\gamma)-z^m\text{Id}_{\hi} = \sum_{k=0}^{m-1}z^k\met^{m-1-k}(\gamma)(\met(\gamma)-z\text{Id}_{\hi})$ and since $\met(\gamma)$ is a unitary operator,
\[
\norm{(\met^m(\gamma)-z^m\text{Id}_{\hi})u}{}{} \le m\norm{(\met(\gamma)-z)u}{}{}
\]
}
and  taking a $z$ minimizing $\norm{(\met(\gamma)-z\text{Id}_{\hi})u}{\hi}{}$ proves the proposition.
\end{proof}
Combining \autoref{esta1} and \autoref{estA} gives
\begin{lem}\label{prop:A1j}
There exists $C_1'\in\RR$ such that for every $m\in\ZZ$,
\[
\norm{A_{1,N}(m) u}{}{} \le C_1'\norm{A_N u}{}{} +|m|\norm{A_{1,N}}{\cal{B}(\hi)}{}\min_{z\in\CC:|z|=1}\norm{(\mathcal{M}_N(\gamma)-z)u}{\hilb_{N}}{}+O(N^{-\8})\norm{u}{\hi}{}.
\]
\end{lem}
{{$A_{\igrk,N}$ corresponds to a long logarithmic time (precisely $2T'$ is asymptotically twice the Ehrenfest time), thus the symbols $a_w$ are no longer in a comfortable symbol class. To solve this issue we first introduce estimates on the norm of $A_{\mathcal Z,N}{}{}$ which is the sum of operators corresponding to the shorter time $T$, for which the pseudodifferential calculus still holds.} We recall the almost-monotonicity property from {\cite[lemma 4.5]{DJ17} which holds for operators $A_{c,N}:=(A_c)_N$,  
\begin{lem}[Almost-monotonicity property]\label{lemono}
	Assume that $c_1,c_2:\mathcal{W}(T)\to\RR$ and for every $w\in \mathcal{W}(T)$, $\abs{c_1(w)}\le c_2(w)\le 1$ .Then there exists a constant $C\in\RR$ {independent of $N, c_1$ and $c_2$} satisfying for every $N\ge 1$ and every $u\in\hi$, 
	\[
	\norm{A_{c_1,N}u}{\hi}{} \le \norm{A_{c_2,N}u}{\hi}{} + CN^{-\frac18}\norm{u}{\hi}{} ,
	\]
	with $C\in\RR$ independent of $c_1,c_2$.
\end{lem}

\begin{proof}
{The proof follows along the lines of lemma 4.5 in \cite{DJ17}. First, note that plugging \autoref{mineig} into \autoref{lemac}, we can write, restricting to $\hi$, 
\[
\Op_N(a_{c_1}) = A_{c_1,N} + O(h^{\frac 14}) \qquad \Op_N(a_{c_2}) = A_{c_2,N} + O(h^{\frac 14}).
\]
Thus it is sufficient to prove that 
\[
\norm{\Op_N(a_{c_1})u}{\hi}{2} \le \norm{\Op_N(a_{c_2})u}{\hi}{2} + CN^{-\frac 14}\norm{u}{\hi}{2}
\]

which is equivalent to \[
\inp{\mathscr{A}u,u}_{\hi} \ge -C N^{-\frac 14}\norm{u}{}{2}
\]
with $$\mathscr{A}=\Op_N(a_{c_2})^*\Op_N(a_{c_2})-\Op_N(a_{c_1})^*\Op_N(a_{c_1}).$$ From \autoref{moyalrhorhop},  $\mathscr{A}=\Op_N(\abs{a_{c_2}}^2-\abs{a_{c_1}}^2) + O(N^{-\frac 14})$. Since $\abs{c_1(w)} \le c_2(w)$ we have $0\le \abs{a_{c_2}}^2-\abs{a_{c_1}}^2 \in S_{L^s,\ha,\frac 14}$. Applying G\r arding inequality from \autoref{lemga},
\[
\inp{\Op_N(\abs{a_{c_2}}^2-\abs{a_{c_1}}^2)u,u} \ge -CN^{-\frac 14}\norm{u}{}{}
\]
from which we deduce the lemma.}
\end{proof}

Using the property we prove}
\begin{prop}
\label{propZ} There exists constants $C_1',C_2',C_3'\in\RR$, independent of $\dl$, such that for every $u\in\hi$,
\eq{
\norm{A_{\cal Z,N}u}{\hi}{} \le \frac{C_1'}{\dl}\norm{\Op_{N}(a)u}{}{}+\frac{C_2'\log N}{\dl}\min_{z\in\CC:|z|=1}\norm{(\mathcal{M}_N(\gamma)-z)u}{\hilb_{N}}{}+\frac{C'_3}{N^{\frac18}\dl}\norm{u}{\hi}{}.
}{AZ}
\end{prop}

\begin{proof}
The indicator function $\ind_{\mathcal{Z}}$ satisfies $0\le \dl\ind_{\mathcal{Z}} \le F \le 1,$ with $F$ defined in \eqref{F}, hence from the almost-monotonicity property in \autoref{lemono},
\eq{
\dl \norm{A_{\mathcal{Z},N}u}{\hi}{} \le \norm{A_{F,N}u}{\hi}{}+O(N^{-\frac 18})\norm{u}{\hi}{}.
}{mona}
 
We note that
\eq
{A_{{F,N}} = {\frac{1}{T}} \sum_{j=0}^{T-1}\sum_{w\in \mathcal{W}(T),w_j=1} A_{w,N} = \frac{1}{T}\sum_{j=0}^{T-1} A_{1,N}(j),
}{AF}

From \eqref{mona} and \eqref{AF},  $$\dl\norm{A_{\mathcal{Z},N}u}{}{} \le \norm{A_{F,N}u}{}{} +O(N^{-\frac 18})\norm{u}{}{} \le \frac{1}{T}\sum_{j=0}^{T-1}\norm{A_{1,N}(j)u}{}{} +O(N^{-\frac 18})\norm{u}{\hilb_{N}}{}$$ which by \autoref{prop:A1j} implies \eqref{AZ}.

\end{proof}

{To finish the proof of \autoref{propigrk} we connect between short logarithmic times and long logarithmic ones.}

\begin{proof}[Proof of \autoref{propigrk}]
We observe that $\igrk=\bigsqcup_{j=1}^8 \igrk_j$ with $$\igrk_j=\left\{w=w^{(1)}\cdots w^{(8)}\in \mathcal{W}(8T): \begin{matrix} w^{(1)},\dots,w^{(j-1)}\in \mathcal{W}(T) \\ w^{(j)}\in\mathcal{Z}\\ w^{(j+1)},\dots,w^{(8)}\in\mathcal{W}(T)\setminus \mathcal{Z} \end{matrix} \right\}.$$ {Since $\sum_{w\in \mathcal{W}(n)}A_{w,N}=\text{Id}_{\hi}$ for every $1\le n\le 8T$
\begin{align*}
A_{\igrk_j,N} &= \sum_{w=w^{(1)}\cdots w^{(8)}\in \igrk_j} A_{w^{(1)}\cdots w^{(8)},N} \\
&= \sum_{w\in \igrk_j}A_{w^{(8)},N}(7T)\cdots A_{w^{(j+1)},N}((j+1)T) A_{w^{(j)},N}(jT)A_{w^{(1)}\cdots w^{(j-1)},N} \\ 
&= A_{\mathcal{W}(T)\setminus \mathcal{Z},N}(7T)\cdots A_{\mathcal{W}(T)\setminus \mathcal{Z},N}((j+1)T) A_{\mathcal{Z},N}(jT)\sum_{w^{(1)},\cdots, w^{(j-1)}\in \mathcal{W}(T)} A_{w^{(1)}\cdots w^{(j-1)},N}\\ 
&=A_{\mathcal{W}(T)\setminus \mathcal{Z},N}(7T)\cdots A_{\mathcal{W}(T)\setminus \mathcal{Z},N}((j+1)T)A_{\mathcal{Z},N}(jT).
\end{align*} 
} {We note that applying \autoref{lemono} there exist constants $C,C'>0$,independent of $N$, such that
\[
{\norm{A_{\mathcal{W}(T)\setminus\mathcal{Z},N}}{}{} \le} \norm{\text{Id}_{\hi}}{\hi}{}+C'N^{-\frac 18} \le C
\]
and therefore for every function $u\in\hi$,
\[
\norm{A_{\igrk,N}u}{\hi}{} \le C\sum_{j=1}^8\norm{A_{\igrk_j,N}u}{\hi}{} \le 8C\norm{A_{\cal Z,N}u}{\hi}{}.
\]
}From \autoref{propZ}, there are constants $C_1,C_2\in\RR$ such that
\[
\norm{A_{\igrk,N}u}{}{} \le \frac{C_1}{\dl\log|\lm_u|}\norm{\Op_{N}(a)u}{}{}+16C\frac{\log N}{\dl\log|\lm_u|}\min_{z\in\CC:|z|=1}\norm{(\mathcal{M}_N(\gamma)-z)u}{\hilb_{N}}{}+\frac{C_2}{N^{\frac18}}\norm{u}{\hi}{}.
\]
\end{proof}

\section{$\nu$-porous sets used for proving \autoref{propfup}}\label{fups}
The next two sections are dedicated to deducing \autoref{propfup} from a version of the fractal uncertainty principle presented in \cite{DJ17}. We first prove a similar result for the Weyl quantization $\Op_h$ acting on $L^2(\RR)$  and then pass to $\Op_N$ by \ref{lemnksup}.
\subsection{Introducing a partition of unity by smooth cut-offs on $\RR^2$}\label{pty}
We equip $\RR^2$ with the $(y,\eta)$ coordinates. We recall that $\iota$ is a change of coordinates from the $(y,\eta)$-coordinates to $(x,\xi)-$coordinates. {In order to simplify the proofs we consider "twisted" symbols $b_{\bt,t}$ in \eqref{bbtt} supported on {$\iota^{-1}\supp(a_\bt\circ\gamma^t)$ for $t\in\llbracket -T'+1,T'\rrbracket$.}} We construct a smooth partition of unity on $\RR^2$: We tile $\RR^2$ by fundamental cells of the lattice $\iota^{-1}\ZZ^2$: fix a fundamental cell of the lattice $\iota^{-1}\ZZ^2$ and denote it by $S^0$. One can consider its $\ka-$thickeing, denoted as $S^0(\ka)=S^0+\mathbb{B}(0,\ka)$ (where $\mathbb{B}(0,\ka)=\{(y,\eta): y^2+\eta^2<\ka^2\}$).

Denote the projections of $S^0(\ka)$ to the coordinate axes by $\Pi_y(S^0(\ka)),\ \Pi_\eta(S^0(\ka))$ and the diameters of the projections by
 \[
 \ell_y=\max_{y,y'\in \Pi_y S^0(\ka)}|y-y'|,\qquad \ell_\eta = \max_{\eta,\eta'\in\Pi_\eta S^0(\ka)} |\eta-\eta'|.
 \]
We can construct a smooth partition of unity: choose a smooth indicator function $\ind_{S^0}^\ka$ supported inside $S^0(\ka)$, attaining the value 1 on some non-empty subset of $S^0$ and $\sum_{m\in\ZZ^2}\ind_{S^m}^\ka=1$ with $$\ind_{S^m}^{\ka}(y,\eta)=\ind_{S^0}^\ka((y,\eta)-\iota^{-1}m).$$ Then the periodization of $\ind_{S^0}^\ka$ gives such partition.
For every  $a\in C^\8(\RR^2)$ we define its truncation on $S^m(\ka)=S^0(\ka)+\iota^{-1}m$ by
\eq{
{}_{m}a = a\cdot \ind^{\ka}_{S^{m}}.
}{ma}

 We recall the notion of $\nu-$porous subsets of $\RR$:
\begin{defn}
	Let $\nu\in(0,1)$ and $0\le \tau_0\le \tau_1$. A subset $\Omg\sub\RR$ is $\nu-$porous on scales $[\tau_0,\tau_1]$ if for each interval $I$ of Lebesgue measure $|I|\in[\tau_0,\tau_1]$ there exists a sub-interval $J\sub I$ such that $|J|=\nu|I|$ and $J\cap\Omg=\emptyset$.
\end{defn}

Recalling \eqref{pra} and $\bt\in\{1,2\}$ we denote the twist of the symbol $a_\bt$ evolved at time $t$ by 
\eq{
b_{\bt,t} = a_\bt\circ\gamma^{t}\circ\iota.
}{bbtt}

Let us recall that $T'=4T$ (where $T$ is defined in \eqref{T}) and consider long words $w\in \mathcal{W}(2T')$. Every such word can be written as concatenation of two words of length $T'$, that is $w=w^-w^+$ with $w^\pm\in\mathcal{W}(T')$,
\begin{align}\label{bpm2}
w^+ &= w_{T'}^+\dots w_{1}^+\quad  w^- = w_{0}^-\dots w^-_{-T'+1}\quad w_{k}^\pm=w_{k+T'-1}\quad \begin{matrix*}
-T'+1&\le {k}&\le T'
\end{matrix*},
\end{align}
and the associated symbols
\begin{align}
b_{w}^+ &= \prod_{{k}=1}^{T'}b_{w_{k}^+,{k}} = \prod_{{k}=1}^{T'} a_{w_{{k}}^+}\circ\gamma^{{k}}\circ\iota,\qquad b_{w}^- =\prod_{{k}=-T'+1}^{0} b_{w_{{k}}^-,{k}} =  \prod_{{k}=-T'+1}^{0} a_{w_{{k}}^-}\circ\gamma^{{k}}\circ\iota.\nonumber
\end{align}

{ We recall that ${}_m b_w^\pm$ is the truncation of $b_w^\pm$ in $S^m(\kappa)$ (see \eqref{ma}), and we define the associated projections by
\[
{}_m\Omg_{w}^{+} = \Pi_y(\supp {}_mb_{w}^+),\qquad {}_m\Omg_{w}^{-} = \Pi_\eta(\supp {}_mb_{w}^-).
\]
Recalling that the $y-$axis corresponds to the unstable direction of $\iota^{-1}\circ\gamma\circ\iota$ and that the $\eta-$axis corresponds to the stable direction of it we prove the following central lemma: 

\begin{prop}\label{concr}
	There exist $K>0$ and $\nu\in (0,1)$ such that the sets ${}_m\Omg_{w}^\pm$ are $\nu-$porous on scales $[Kh^\rho,1]$.
\end{prop}

\begin{proof}
	We prove the porosity statement only for ${}_m\Omg_{w}^+$ as the proof for ${}_m\Omg_{w}^-$ follows along the same lines. For simplicity we take $m=0$ as the proof for the general case is analogous and follows by translation. {The proof relies on the correspondence between finite intervals in $\RR$ and unstable sub-orbits of $\gamma$.}
	\begin{enumerate}[label=\arabic*., leftmargin=*, ref=part \arabic*]
	
	\item { For every fixed $(v_0,t_0)\in\TT^2\times\RR$ we define the truncated unstable orbit of length $\tau$ by 
	\eq{
	\mathcal{O}_{v_0,t_0,\tau} := \{\mathcal{H}_u(v_0,t): t\in[t_0,t_0+\tau]\}.
}{trn_orbt}
We denote for brevity $\mathcal{O}_{v_0,\tau}:=\mathcal{O}_{v_0,0,\tau}$. Denote the ball of radius $r$ centered at $(y,\eta)$ by $\mathbb{B}((y,\eta),r)$. Let us recall from the construction above that (cf. \autoref{fig:pt}) \eq{\mathcal{K}_1\sub\TT^2\setminus\supp(a_1),\qquad \mathcal{K}_2\sub\TT^2\setminus\supp(a_2)}{Ksets}  are both open hence there exists $\ka'>0$ and balls $\mathbb B_1((y_1,\eta_1),\ka') , \mathbb B_2((y_2,\eta_2),\ka')\in\TT^2$ such that $$ \mathbb B_1((y_1,\eta_1),\ka')\sub\mathcal{K}_1,\qquad \mathbb B_2((y_2,\eta_2),\ka')\sub\mathcal{K}_2.$$ We note that each ball $\mathbb B_\bt((y_\bt,\eta_\bt),\ka')$ contains a smaller ball $\mathcal{K}'_\bt:=\mathbb B((y_\bt,\eta_\bt),\frac{\ka'}{4})$.
Henceforth let us denote the Euclidean distance (both in $\RR^2$ and in $\TT^2$) by $\text{d}_{Euc}$. From the construction
\eq{
\text{d}_{Euc}(\bnd \mathcal{K}_\bt,\mathcal{K}'_\bt) \ge \frac 34 \ka'.
}{d>34}
} 
 	\item\label{sod} { Let us prove that if a truncated unstable orbit is long enough it has to pass through both sets $\mathcal{K}_\bt'$ and  moreover that the lengths of the truncated orbits lying inside $\mathcal{K}_\bt'$ are uniformly bounded from below. In other words, we show that:
 	
There exist some $(L,\ell)\in(\RR_+^*)^2$ with $L>\ell$ such that for every $v_0$ the truncated unstable orbit $\mathcal{O}_{v_0,L}$ contains sub-orbits $\mathcal{O}_{v_0,t_1,\ell},\mathcal{O}_{v_0,t_2,\ell}$ satisfying $\mathcal{O}_{v_0,t_j,\ell}\sub\mathcal{K}_j'$.
 	
\noindent Assume to the contrary that for every $(L,\ell)$ there exists an orbit $\mathcal{O}_{v_0,L}$ which does not contain  either a sub-orbit $\mathcal{O}_{v_0,t_1,\ell}\sub\mathcal{K}_1'$ or a sub-orbit $\mathcal{O}_{v_0,t_2,\ell}\sub\mathcal K_2'$. Let us consider the sequences $L_n=n$ and $\ell_n=\frac 1n$ and construct a sequence of sub-orbits of lengths $L_n$: for every $n$ there exists a vector $v_n\in\TT^2$ such that the orbit $\mathcal{O}_{v_n,n}$ does not contain a sub-orbit $\mathcal{O}_{v_n,t_\bt,\frac 1n}\sub\mathcal{K}_\bt'$. Since $\TT^2$ is compact there exist a convergent sub-sequence of vectors $v_{n_k}\to v_{\8}\in\TT^2$. Applying the assumption, for every $k\in\NN$ there exists $\bt_k\in\{1,2\}$ such that the orbit $\mathcal{O}_{v_{n_k},n_k}$ does not contain a sub-orbit $\mathcal{O}_{v_{n_k},t_{\epsilon_k},n_k^{-1}}\sub\mathcal{K}_{\epsilon_k}'$. Possibly extracting a sub-sequence, we can ensure that for every $k\ $ $\bt_k=\bt$. Letting $k\to\8$ the sub-orbits $\{\mathcal{O}_{v_{n_{k}},n_{k}}\}_k$ converge to a full orbit $\mathcal{O}_{v_{\8}}$ which does not intersect $\mathcal{K}_\epsilon'$. Since the latter is a proper open subset of $\TT^2$, it is a contradiction to} the minimality of the horocylic flow we proved in \autoref{dens} and thus we have proved the claim. 
 	
 	\item Fix some big constant $K>0$ and let $I=I(h)\sub\RR$ be an interval with $\abs{I} \in [Kh^{\rho},1]$ and let \eq{k_0 =\left\lceil \frac{\log (8\sqrt{(m_u-m_s)(1+m_s^2)}(1+\ka))-\log(3\ka')}{\log\abs{\lm_u}}\right\rceil.}{k0} Recall \eqref{T} and define 
 	\eq{
 		k = \min_{k_0< j< T'}\{j:\abs{\lm_u}^j\abs{I} \ge L\}.
 	}{k}

\item\label{soq} {Let us associate to $I$ a truncated unstable orbit:} Recalling the projection $\pi$ from \eqref{pi} we consider the set 
 \eq{\mathcal{O}_I:= \left\{\pi\circ\iota\binom y0:y\in I\right\} = \{\pi(y.v_u):y\in I\} \sub\TT^2. 
}{oi}
From {\autoref{irr}} it is indeed a truncated unstable orbit of length $\abs{I}$. 

{
\item Recalling \eqref{Hin} implies that $\gamma^k\mathcal{O}_I$ is an orbit of length $\abs{\lm_u}^k\abs{I}\ge L$. 
 From \ref{sod} and \ref{soq} there exists a sub-orbit $\mathcal{O}'_I\sub \gamma^k\mathcal{O}_I\cap \mathcal{K}_{w_k^+}'$ { of length $\ell$. Then there exists a sub-interval $J\sub I$ such that $\gamma^{-k}\mathcal{O}_I' = \{\pi\circ\iota\binom y0:y\in J\}$ and $|J|=\abs{\lm_u}^{-k}\ell$.
\item For establishing the porosity let us first show that there exists some $\nu\in(0,1)$} such that $\frac{\abs{J}}{\abs{I}}\ge \nu$.
{ If $k > k_0+1$, using the minimality of $k$
\[
\abs{I}\abs{\lm_u}^k < \abs{\lm_u} L,
\] which implies
\eq{
\frac{\ell}{L}\ge\frac{\abs{J}}{\abs{I}} = \frac{\ell\abs{\lm_u}^{-k}}{\abs{I}}\ge \frac{\ell}{\abs{\lm_u}L}.
}{k>k0+1}
Similarly if $k=k_0+1$, 
\eq{
\frac{\ell}{L}\ge\frac{\abs{J}}{\abs{I}} = \frac{\ell\abs{\lm_u}^{-k_0-1}}{\abs{I}}\ge {\frac{\ell}{\abs{\lm_u}^{k_0+1}}}.
}{k=k0+1}
As follows, taking $\nu =\frac{\ell}{\abs{\lm_u}^{k_0+1}L}$ both \eqref{k>k0+1} and \eqref{k=k0+1} are satisfied.
} }
\item {For demonstrating the $\nu-$porosity of ${}_0\Omg_{w}^+$ let us prove that for the choice of $k_0$ in \eqref{k0} $J$ does not intersect ${}_0\Omg_{w}^+$. } 
Assume to the contrary there is some $y\in J\cap{}_0\Omg_{w}^+$ then it can be lifted to a point $(y,\eta)\in\supp {}_0 b_w^+$ for some $\eta\in\RR$. Let us write, 
\eq{
\gamma^k\circ\pi\circ\iota\binom y\eta = \gamma^k\circ\pi\circ\iota\binom y0+\gamma^k\circ\pi\circ\iota\binom 0\eta.
}{dec}
From the construction of $J$ we deduce that the first summand is in $\mathcal{K}_{w_k^+}'$. Since $(y,\eta)\in S^0(\ka)$ we can bound $|\eta|\le \ell_\eta$. Let us bound the second summand as well: { We recall that $\iota^{-1}$ maps any truncated stable orbit $\mathcal{O}_{v_0,t_0,\tau}^s:=\{\mathcal{H}_s(v,t):t\in I'=[t_0,t_0+\tau]\}$ to an interval $I''\sub \{c\}\times\RR$ (for some $c\in\RR$) and that $$d_{\text{Euc}}(\iota^{-1}\mathcal{H}_s(v_0,t_0),\iota^{-1}\mathcal{H}_s(v_0,t_0+\tau))=\norm{\iota^{-1}v_s}{}{}\tau.$$ We note that $$\abs{\left\{\gamma^k\circ\pi\circ\iota\binom 0\eta:|\eta|\le \ell_\eta\right\}}=2\ell_\eta\abs{\lm_u}^{-k}\norm{v_s}{}{}$$ and
\[
\abs{\iota^{-1}\circ\gamma^k\circ\iota\binom 0\eta} \le 2\frac{\ell_\eta}{\sqrt{m_u-m_s}}\abs{\iota^{-1}\circ\gamma^k v_s}  = 2\ell_\eta\abs{\lm_u}^{-k}\norm{v_s}{}{}.
\]
\begin{figure}[tp]
\centering
\includegraphics[scale=0.4]{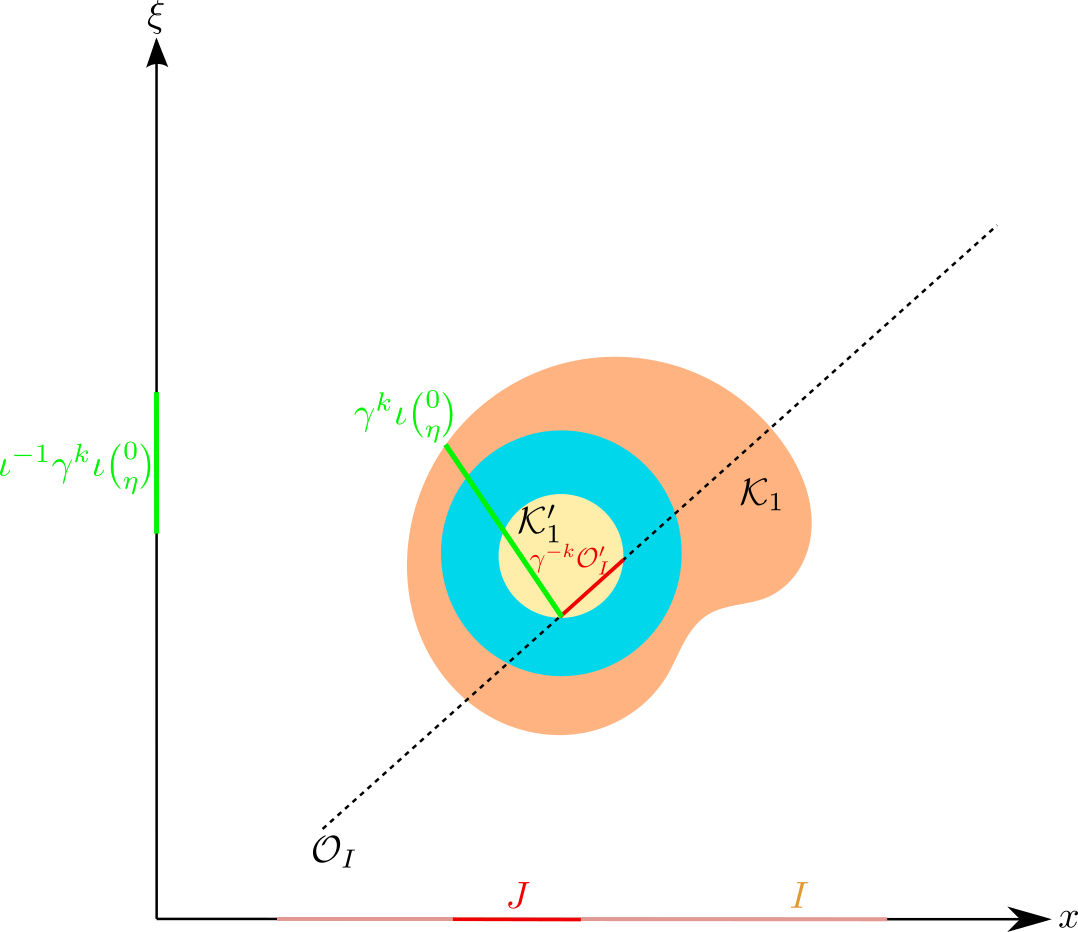}
\caption{The ball of radius smaller than $\frac 34\ka'$ centered at $\gamma^k\circ\pi\circ\iota\binom y\eta$ from \eqref{dec}, appearing in the proof of \autoref{concr}, is contained inside $\mathcal{K}_{w_k^+}$. In the figure $w_k^+=1$.}
\label{figdes2}
\end{figure}
As follows the length of the truncated stable orbit $\{\gamma^k\circ\pi\circ\iota\binom 0\eta:\abs{\eta}\le \ell_\eta\}$ is bounded from above by $2\ell_\eta\abs{\lm_u}^{-k_0}\norm{v_s}{}{}$. } From the choice of $k_0$ in \eqref{k0} the second summand in \eqref{dec} is bounded from above by $\frac{3}{4}\ka'$ and from \eqref{d>34} (see \autoref{figdes2}) adding it to the first summand, $\gamma^k\circ\pi\circ\iota\binom y\eta\in \mathcal{K}_{w_k^+}\sub\complement \supp(a_{w_k^+})$ or equivalently $(y,\eta)\notin \supp(b_{w}^+)$ which is a contradiction. 
\end{enumerate}
	\end{proof}
}
\section{An estimate on words in the uncontrolled region}\label{finalfups}
We begin this section by connecting between the operator $A_{w}$ and $\Op_h(b_w^\pm)$: The quantization of $b_w^\pm$ satisfies a central role in estimating $\norm{A_w}{\bltr}{}$ for $w\in\mathcal{W}(2T')$. Since $2T'$ is bigger than Egorov time $\frac{\log\frac 1h}{\log|\lm_u|}$, we cannot view $A_w$ as a pseudo-differential operator having a symbol in the previously defined calculi. {Nevertheless we show that $\norm{A_w}{}{}$ equals to the norm of a product of two pseudo-differential operators up to some reminder}. Recalling \autoref{remprodop}, we will deduce an estimate on $\norm{A_{w,N}}{\mathcal{B}(\hi)}{}\le \norm{A_w}{\bltr}{}$, proving \autoref{propfup}. Using the unitarity of $\mathcal{M}_h(\gamma^{T'-1})$ and \autoref{lemaw} 
\begin{align}\label{Aws}
\norm{A_w}{\bltr}{} &= \norm{A_{w_{2T'-1}}(2T'-1)A_{w_{2T'-2}}(2T'-2)\cdots A_{w_{T'}}(T') \cdots A_{w_1}(1)A_{w_0}}{\bltr}{}\\
&=\bigg\Vert\mathcal{M}_h(\gamma^{T'-1})A_{w_{2T'-1}}(2T'-1)\cdots 
A_{w_1}(1)A_{w_0}\mathcal{M}_h(\gamma^{-T'+1})\bigg\Vert_{\bltr}\nonumber \\
&=\bigg\Vert A_{w_{2T'-1}}(T')\cdots A_{w_{T'}}(1) \cdots A_{w_1}(-T'+2)A_{w_0}(-T'+1)\bigg\Vert_{\bltr}\nonumber\\
&= \bigg\Vert A_{w_{T'}^+}(T')\cdots A_{w_{1}^+}(1)A_{w_0^-}(0) \cdots A_{w_{-T '+2}^-}(-T'+2)A_{w_{-T'+1}^-}(-T'+1)\bigg\Vert_{\bltr}\nonumber
\end{align}

Invoking the unitary of $\mathcal{M}_h(\iota)$ and \autoref{lemaw}
\begin{align}\label{prodb}
&= \bigg\Vert \mathcal{M}_h(\iota)^{-1}A_{w_{T'}^+}(T')\cdots A_{w_1^+}(1) \cdots A_{w_{-T'+2}^-}(-T'+2)A_{w_{-T'+1}^-}(-T'+1)\mathcal{M}_h(\iota)\bigg\Vert_{\bltr} \\ 
&=\norm{(\Op_h(b_{w}^+)+O(h^{1-\rho-\eps}))(\Op_h( b_{w}^-)+O(h^{1-\rho-\eps}))}{\bltr}{}\nonumber\\
&\le \norm{\Op_h(b_{w}^+)\Op_h(b_{w}^-)}{\bltr}{} + O(h^{1-\rho-\eps}),\nonumber 
\end{align}

\subsection{Proof of \autoref{propfup}} In this subsection we combine the porousity of the sets ${}_m\Omg_{w}^{+},{}_m\Omg_{w}^{-}$ proved in \autoref{concr} with the fractal uncertainty principle in order to deduce \autoref{propfup}. For every $m,m'\in\ZZ^2$ let us consider
\[
B_{(m,m')}=\Op_h({}_{m}b_{w}^+)\Op_h({}_{m'}b_{w}^-).\]	
We will apply the Cotlar-Stein theorem (cf. \autoref{stein} below) in order to prove that
\begin{prop}\label{global}
{There exists $\be'>0$ such that for every two words $w^\pm\in\mathcal{W}(T')$} 
\[
\norm{\Op_h(b_{w}^+)\Op_h(b_{w}^-)}{\mathcal{B}(\Ltr)}{} = O(h^{\be'}).
\]
\end{prop} 
Combining \autoref{global},\eqref{prodb}, the unitarity of $\mathcal{M}_h(\iota)$ and \ref{lemnksup} we deduce \autoref{propfup} for every $w\in\mathcal{X}$
	\[
	{
			\norm{A_{w,N}}{\mathcal{B}(\hi)}{} \le \norm{A_w}{\bltr}{}  \le Ch^{\be'}. 
	}
	\]
where $\be'=\min\{\be,1-\rho\}$ for $\be$ defined in \autoref{propfup2} below.
For the rest of the section we prove \autoref{global}. Let us denote
\eq{
\ell(\ka)=\max\{\ell_y,\ell_\eta\}.
}{ell}
  We begin by obtaining local estimates. Let us define for every $m,m'\in\ZZ^2$ $${}_m{B_-}= \Op_h({}_mb_w^-),\ {}_m{B_+}=\Op_h({}_mb_{w}^+).$$ We will provide two distinct estimates on $\norm{B_{(m,m')}}{\Ltr}{}$ based on the distance between $m$ and $m'.$ First we show that when $m$ is in the vicinity $m
'$ this norm is $O(h^\be)$.

\subsubsection{Local estimates corresponding to nearby fundamental cells}
We recall the concept of smooth cutoffs on a compact $I\sub\RR$. Consider the $K-$thickened set $I(K)=I+[-{K},{K}]$. {We consider the smooth cut-off $\chi_{I}^{{K}}(x)\in C_c^\8(\RR)$ supported on $I(2 K)$, equal to unity in $I(K)$ and $\chi_I^{{K}}\mid_{(I(2 K))\setminus I( K)}(x)\in \left[0,1\right]$.} The non-smooth characteristic function of a compact $I \sub\RR$ will be denoted by $\ind_I$.

For the class of $\nu$-porous sets in $\RR$ Dyatlov and Jin have established in \cite{DJ17} (based on preceding \cite{BD16}) an extension of the "original" uncertainty principle. We will use the following version of it, 
\begin{prop}[Fractal uncertainty principle, {\cite[Proposition 2.10]{DJS22}}]
	\label{propfup2}
	Fix some $K>0, \rho \in (\ha,1)$ and $\nu\in(0,1)$. Then there exist $\be>0$ depending on $\nu,\rho$ and a constant $C$ depending on $\nu,K$ such that for every $\nu-$porous sets $X,Y\sub \RR$ on scales $Kh^\rho$ to $1$,
	\[
	\norm{\ind_{Y(2Kh^\rho)}\F_h^{-1} \ind_{X(2Kh^\rho)}}{\Ltr\to\Ltr}{} \le Ch^\be\qquad h\in(0,1),
	\] 
	where $X(2Kh^\rho)=X+[-2Kh^\rho,2Kh^\rho],Y(2Kh^\rho)=Y+[-2Kh^\rho,2Kh^\rho]$.
\end{prop}
The connection between $\be$ to $\rho$ is expressed in \cite{DJS22} as $\be=\be_0(2\rho-1)$ for some constant $\be_0>0.$
\begin{lem}\label{larem}
	The sharp cutoffs $\ind_{Y(2Kh^\rho)}$ and $\ind_{X(2Kh^\rho)}$ in \autoref{propfup2} can be replaced by their smooth counterparts $\chi_Y^{Kh^\rho}$ and $\chi_X^{Kh^\rho}$. Indeed,
	\begin{align*}
	\norm{\chi_{Y}^{Kh^\rho}\F_h^{-1} \chi_{X}^{Kh^\rho}}{\Ltr\to\Ltr}{} &= \norm{\chi_{Y}^{Kh^\rho}\ind_{Y(2Kh^\rho)}\F_h^{-1} \ind_{X(2Kh^\rho)}\chi_{X}^{Kh^\rho}}{\Ltr\to\Ltr}{}  \\&\le \norm{\ind_{Y(2Kh^\rho)}\F_h^{-1} \ind_{X(2Kh^\rho)}}{\Ltr\to\Ltr}{} \le Ch^\be
	\end{align*}
\end{lem}

Suppose the intersection of the fundamental cells is non-empty we deduce a bound on the norm $\norm{{}_mB_{+}{}_{m'}{B_{-}}}{}{}$ in virtue of the fractal uncertainty principle, or rather its formulation in \autoref{larem} above,
\begin{lem}\label{portoy} 
Let $m,m'\in\ZZ^2$ then there exists $\be>0$, independent of $m$ and $m'$, and a constant $C'>0$ such that for every $h\in\left(0,1\right]$
	\[
	\sup_{w=w_+w_-\in\mathcal{W}(2T')}\norm{{}_mB_{+}\ {}_{m'}{B_{-}}}{}{} \le C'h^{\be}
	\] 
	with $T'$ specified in \eqref{T}, ${}_m B_+=\Op_h({}_mb_w^+)$ and ${}_{m'}B_-=\Op_h({}_{m'}b_w^-)$, the operators corresponding to halves of a word $w\in\mathcal{W}(2T')$ (see \eqref{bpm2}).
\end{lem}

\begin{proof}
Recall that $\chi^{Kh^\rho}_{{}_m\Omg_{w}^{+}}\in C_c^\8(\RR)$ is unity on ${}_m\Omg_{w}^{+}(Kh^\rho)$ and that $\Pi_y\supp({}_mb_+)\sub {}_m\Omg_{w}^{+}$. Then since $\chi^{Kh^\rho}_{{}_m\Omg_{w}^{+}}=\Op_h(\chi^{Kh^\rho}_{{}_m\Omg_{w}^{+}})$ and {since we can take $\chi^{Kh^\rho}_{{}_m\Omg_{w}^{+}}\ \in S_{L^\xi,\rho,0}(\RR^2)$}, from the Moyal product in the class $S_{L^\xi,\rho,0}$,
\[
{}_{m} B_+\chi^{Kh^\rho}_{{}_m\Omg_{w}^{+}}  = {}_{m} B_+ + \Op_h(r_+).
\]

In fact, since all terms of the Moyal product of ${}_{m} B_+\chi^{Kh^\rho}_{{}_m\Omg_{w}^{+}}$  except the first vanish repeated integration by parts as in \autoref{nsp} implies that the norm of the reminder is $\norm{\Op_h(r_+)}{}{}=O(h^\8)$. Analogously, since $\chi_{{}_{m'}\Omg_-}^{Kh^\rho}\in C_c^\8(\RR)$ and $\Pi_\eta\supp({}_{m'}b_-)={}_{m'}\Omg_-$
\[
\F_h^{-1}\chi^{Kh^\rho}_{{}_m\Omg_-}\F_h\ {}_{m'} B_-  = {}_{m'} B_- + O_{\bltr}(h^\8).
\]

As a result, 
\begin{align*}
	\norm{{}_m B_+\ {}_{m'} B_-}{}{} &\le \norm{{}_m{B_+}\chi_{{}_{m}\Omg_{w^+}}^{Kh^\rho} 
		\F_h^{-1}\chi_{{}_{m'}\Omg_{w^-}}^{Kh^\rho}\F_h\ {}_{m'}{B_-}}{}{}+O(h^{\8}) \\
	&\le \norm{{}_mB_+}{}{}\norm{{}_{m'}B_-}{}{}\cdot \norm{\chi_{{}_{m}\Omg_{w^+}}^{Kh^\rho} 		\F_h^{-1}\chi_{{}_{m'}\Omg_{w^-}}^{Kh^\rho}}{}{}+O(h^{\8}).
	\end{align*}
There exists $\nu\in (0,1)$ and a large $K$ such that ${}_m\Omg_{w}^{+}$ and ${}_{m'}\Omg_{w}^-$ are $\nu-$porous sets on scales $[Kh^\rho,1]$.  {From \autoref{bdopa} the norms of ${}_mB_\pm$ are uniformly bounded}. Applying \autoref{larem} there is some $\be>0$ such that 
\[
\norm{{}_mB_{+}\ {}_{m'}B_{-}}{}{} \le C'' \norm{\chi_{{}_m\Omg_{w}^{+}}^{Kh^\rho}\F_h^{-1}\chi_{{}_{m'}\Omg_{w}^-}^{Kh^\rho}}{}{}+O(h^{\8}) \le C'h^{\be}.
\]
\end{proof}

\subsubsection{Local estimate for distant fundamental cells}
We next provide a local estimate when the two cutoffs are disjoint: 
\begin{lem}\label{lemhinty}
	Let $k\in\NN^*$ be a large integer, $\rho\in \left[0,1\right)$ and $m,m'\in\ZZ^2$ such that $\norm{\iota^{-1}(m-m')}{}{}\ge 10\ell(\ka)$. Then 
	
	there exist a universal constant $M\in\NN$ and a scalar $C_k$  independant of $m, m'$ such that for every $\mathscr{a,b}\in C^\8(\RR^2)$ satisfying $\mathscr{a},\mathscr{b} \in S_\rho(\RR^2)$ 
	\[
	\norm{\Op_h({}_m{\mathscr{a}})\Op_h({}_{m'}{\mathscr{b}})}{}{} \le C_k\frac{N_{k+M}(\mathscr a)N_{k+M}(\mathscr b)}{\norm{m-m'}{}{k}}h^{k(1-\rho)+M(\ha-\rho)-2},
	\]
	where $N_k$ is defined by
	\[
	N_k(\mathscr a) :=\max_{\abs{\al}\le k}\sup_{h\in \left(0,1\right]}\left(h^{\rho\abs{\al}}\norm{\bnd^\al \mathscr a}{\8}{}\right).
	\]
\end{lem}

The proof of the lemma is a direct adaptation of \cite[Theorem 4.11.ii]{Z12} followed by  a Calderon-Vaillancourt estimate,
\begin{proof}[Sketch of proof]
	\begin{enumerate}[label=\arabic*., leftmargin=*, ref=part \arabic*]
\item[]
\item Let us denote the symplectic form on $\RR^2$ by $\sigma$.  
Denote by ${}_{m,m'}{\mathscr c}$ the Moyal product of the symbols, i.e., 
\begin{align*}\label{WeylF}
	\Op_h({}_{m,m'}{\mathscr c})=\Op_h({}_m{\mathscr{a}})\Op_h({}_{m'}{\mathscr{b}}). 
\end{align*}	
Recall the integral representation of ${}_{m,m'}{\mathscr c}$ given by  
\eq{
	{}_{m,m'}{\mathscr c}(Z) = \mathcal{I}(Z) = \frac{1}{(\pi h)^2 }\int_{\RR^2}\int_{\RR^2} e^{-\frac{2i}{h}\sigma(z_1,z_2)}{}_m {\mathscr{a}}(Z-z_1){}_{m'}{\mathscr b}(Z-z_2) dz_1dz_2.
}{bmm'1}
\item Let us bound the derivatives $\bnd^{\gamma}_Z{}_{m,m'}\mathscr c$. Due to the support properties of ${}_m\mathscr a$ and ${}_{m'}\mathscr b$, the integrand in \eqref{bmm'1} is supported on a set of bounded volume away from the origin and as a result one can show that $\norm{z}{2}{} = \norm{(z_1,z_2)}{2}{} \ge 8\ell(\ka)$. We deduce that $\Vert\nabla\sigma\Vert$ does not vanish anywhere on its support as
\[
\norm{\nabla\sigma}{2}{} = \norm{(\eta_2,-y_2,-\eta_1,y_1)}{2}{} = \norm{z}{2}{} \ge 8\ell(\ka).
\]
 Define now the differential operator \eq{\mathscr{D}(a)=-\frac{h}{2i}\frac{\inp{\nabla\sigma,\nabla a}}{\norm{\nabla\sigma}{}{2}},\qquad \mathscr{D}:C^\8(\RR^4)\to C^\8(\RR^4).}{D}
Then since $\nabla\sigma(z)$ is a linear function of $z$, for every $k\in\NN^*,$
\[
\abs{({}^t\mathscr{D}^ka)(z)}\le Ch^k\sum_{j=0}^k\frac{\norm{\bnd_z^{j}a}{}{}}{\norm{z}{}{2k-j}}.
\]
For every $k\in\NN^*$ we recall that since $\mathscr{a},\mathscr{b}\in S_\rho(\RR^2)$, $$\norm{\bnd^j_z ({}_m\mathscr{a}(Z-z_1){}_{m'}\mathscr{b}(Z-z_2))}{}{}\le C_jh^{-\rho j}$$ thus
\eq{
	\abs{{}^t(\mathscr{D}^k)({}_m\mathscr{a}(Z-z_1)\ {}_{m'}\mathscr{b}(Z-z_2))}\le C h^k (\sum_{j=0}^k\frac{\norm{\bnd^{j}({}_{m}\mathscr{a}\ {}_{m'}\mathscr{b})}{}{}}{\norm{z}{}{2k-j}}) \le C_k\frac{h^{k(1-\rho)}}{\norm{z}{}{k}}N_k(\mathscr a)N_k(\mathscr b).
}{bd1}
\item Let us estimate the integral appearing in \eqref{bmm'1} through integration by parts. We recall that the integral is supported on a product of two compact domains thus we can apply change of variable
\[
z_1 := \iota^{-1}m-Z+c_1, \  z_2 := \iota^{-1}m'-Z+c_2
\]
for some $c_1,c_2=O(1)$. Using the Pythagorean law and then the reverse triangle inequality
\begin{align}
	\abs{{}^t(\mathscr{D}^k)({}_m\mathscr{a}(Z-z_1){}_{m'}\mathscr{b}(Z-z_2))} &\le C_k\frac{h^{k(1-\rho)} N_k(\mathscr a)N_k(\mathscr b)}{(\norm{z_1}{}{2}+\norm{z_2}{}{2})^{\frac k2}} \\&=  C_k\frac{h^{k(1-\rho)}N_k(\mathscr a)N_k(\mathscr b)}{(\norm{\iota^{-1}m-Z+c_1}{}{2}+\norm{\iota^{-1}m'-Z+c_2}{}{2})^{\frac k2}} \nonumber\\
	&\hspace{-6em}\le\frac{C'N_k(\mathscr a)N_k(\mathscr b)}{((\norm{\iota^{-1}(m-m')}{}{}-2\ell(\ka))^2+4(\norm{Z-\ha\iota^{-1}(m+m')}{}{}-2\ell(\ka))^2)^{\frac k2}}h^{k(1-\rho)}\nonumber
\end{align}
and  therefore integraring over $c_1,c_2$ we get
\begin{align*}
	\abs{{}_{m,m'}{\mathscr c}(Z)}&=\abs{\mathcal{I}(Z)} \\&\le  \frac{C''(N_k(\mathscr a), N_k(\mathscr b))\text{Vol}(S^0(\ka))^2}{((\norm{\iota^{-1}(m-m')}{}{}-2\ell(\ka))^2+4(\norm{Z-\ha\iota^{-1}(m+m')}{}{}-2\ell(\ka))^2)^{\frac k2}}h^{k(1-\rho)-2}.
\end{align*}
Finally, in order to apply Calderon-Vaillancourt theorem we establish uniform estimates on a finite number of derivates of $_{m, m'}\mathscr{c}$:
Denoting $C_{m,m',\ell(\ka)} = (\norm{\iota^{-1}(m-m')}{}{}-2\ell(\ka))^2$  and applying integration by parts for every multi-index $\gamma\in\NN^2$,
\eq{
	\abs{\bnd_Z^\gamma {}_{m,m'}{\mathscr c}(Z)}=\abs{\bnd_Z^\gamma \mathcal{I}(Z)} \le\frac{C_k' N_{k+|\gamma|}(\mathscr a)N_{k+|\gamma|}(\mathscr b)}{(C_{m,m',\ell(\ka)}^2+\bigg(\norm{Z-\ha\iota^{-1}(m+m')}{}{}-2\ell\bigg)^2)^{\frac k2}}h^{k(1-\rho)-|\gamma|\rho-2},
}{Igamma}

\item The estimates on the derivatives of ${}_{m,m'}\mathscr{c}$ give an upper bound on $\norm{\Op_h({}_{m,m'}\mathscr{c})}{\bltr}{}$.
We recall from \cite[Theorem 5.1]{Z12} for every $u\in \Ltr$
\[
\Op_h({}_{m,m'}{\mathscr c})u = h^{-\qt}\Op_1({{}_{m,m'}{\check{\mathscr c}}})\check{u},\qquad \check u(\check x)= h^{\frac 14}\check{u}(\sqrt h \check x) , {{}_{m,m'}{\check{\mathscr c}}}(\check y,\check \eta)= {}_{m,m'}{\mathscr c}(\sqrt h\check y,\sqrt h\check\eta)
\]
with the rescaling of $u$ being a unitary operator on $\Ltr$ and $\check c\in S_{\rho-\ha}$. From \cite[Theorem 4.23.ii]{Z12} there exists some universal $M$ such that 
\begin{align*}
	&\norm{\Op_h({}_{m,m'}{\mathscr c})}{\bltr}{} \\&\hspace{4em}=\norm{\Op_1({}_{m,m'}{\check{\mathscr c}})}{\bltr}{} \\
	&\hspace{4em}\le C\sum_{\abs{\gamma}\le M}h^{\frac{\abs{\gamma}}{2}}\norm{\bnd^\gamma {}_{m,m'}{\check{\mathscr c}}}{\8}{} \\
	&\hspace{4em}\le C_k'' \sum_{\abs{\gamma}\le M}\sup_{Z}\frac{ N_{k+|\gamma|}(\mathscr a)N_{k+|\gamma|}(\mathscr b)}{(C_{m,m',\ell(\ka)}^2+\bigg(\norm{\sqrt h Z-\ha\iota^{-1}(m+m')}{}{}-2\ell\bigg)^2)^{\frac k2}}h^{k(1-\rho)+\abs{\gamma}(\ha-\rho)-2}
\end{align*}
The terms whose contributions to the sum are the largest are those corresponding to $\abs{\gamma}=M$. Taking $k$ to be arbitrarily large, we get an $O(h^\infty)$ estimate and

\[
\norm{\Op_h({}_{m,m'}{\mathscr c})}{\bltr}{} \le C_k'''\frac{N_{k+M}(\mathscr a)N_{k+M}(\mathscr b)}{\norm{m-m'}{}{k}}h^{k(1-\rho)+M(\ha-\rho)-2}.
\]
\end{enumerate}
\end{proof}

\begin{rem}\label{remb--}
{Since $S_{L^x,\rho,0}(\RR^2),S_{L^\xi,\rho,0}(\RR^2)\sub S_\rho(\RR^2)$ we deduce that} ${}_m b_w^+,\ {}_{m'}b_w^-\in S_{\rho}(\RR^2)$. For every $m,m'\in\ZZ^2\ $ ${}_m b_w^{\pm},\ {}_{m'}b_w^{\pm}\in S_{\rho}(\RR^2)$ and for every $k\in\NN^*\setminus\{1,2\}$ there exists a constant $C_k$ for which if $\norm{\iota^{-1}(m-m')}{}{}>10\ell(\ka)$
	\begin{align*}
	\norm{{}_mB_\pm\ {}_{m'} B_\mp}{}{} \le h^{k(1-\rho)+M(\ha-\rho)-2}\frac{C_k}{\norm{m-m'}{}{k}},\qquad \norm{{}_mB_\pm^*\ {}_{m'} B_\pm}{}{} \le h^{k(1-\rho)+M(\ha-\rho)-2}\frac{C_k}{\norm{m-m'}{}{k}},
	\end{align*}
	where $M$ is the universal vonstant from \autoref{lemhinty} and ${}_m B_\pm=\Op_h({}_mb_w^\pm),{}_{m'}B_\pm=\Op_h({}_{m'}b_w^\pm)$ are the operators obtained from quantifying the product of $\ind_{S^m}^K$ with a symbol corresponding to a half of a long word $w\in\mathcal{W}(2T')$.
\end{rem}

We deduce now \autoref{global} from the estimates above,
\newline 
\begin{proof}[Proof of \autoref{global}]
Let us begin by recalling Cotlar-Stein theorem
\begin{thm}[Cotlar-Stein theorem, \S 3 in \cite{Co55},\S VII.2.2 in \cite{S97}] \label{stein}
	Let $\{A_j\}_{j\in \NN}$ be a family of bounded operators on some Hilbert space $\hil$. Suppose the  bounds
	\[
	\sup_j\sum_k\norm{A_j^*A_k}{\hil\to \hil}{\ha}\le C,\text{ and }\sup_j\sum_k\norm{A_jA_k^*}{\hil\to \hil}{\ha}\le C 
	\]
	hold then $\sum_j A_j$ converges, in the strong operator topology, to an operator $A$ satisfying $\norm{A}{\hil\to \hil}{}\le C$.
\end{thm}
We will apply this theorem on the family $\{B_{(n,n')}\}_{n,n'\in\ZZ^2}$. From the definition of $b_w^{\pm}$ in \eqref{bpm2} we note that $\norm{\Op_h({}_nb_{\pm})}{}{}\le C$ uniformly in $n\in\ZZ^2$ hence this family is bounded. We would like to show the boundedness of the sums. For that purpose let us define the relation $\sim$ on $\ZZ^2\times \ZZ^2$ by
\[
n\sim n'\iff \norm{\iota^{-1}(n-n')}{}{} \le 10\ell(\ka).
\]
Fix $(q,q')\in\ZZ^4$ and consider the sum $\sum_{n,n'\in\ZZ^2}\norm{B^*_{(q,q')}B_{(n,n')}}{}{\ha}$. Splitting the summation over $n$ and $n'$,
\begin{align*}
\sum_{n,n'\in\ZZ^2}\norm{B_{(q,q')}^*B_{(n,n')}}{}{\ha} &= \sum_{n,n'\in\ZZ^2}\norm{{}_{q'}B^*_-\ {}_{q}B^*_+\ {}_{n}B_+\ {}_{n'}B_-}{}{\ha} \\&= \sum_{n,n':\substack{n\sim n'\\ n\sim q}} \norm{B_{(q,q')}^*B_{(n,n')}}{}{\ha}+\sum_{n,n':n\not\sim n'\text{ or }n\not\sim q}\norm{B_{(q,q')}^*B_{(n,n')}}{}{\ha}.
\end{align*}
Since $q\in\ZZ^2$ is fixed the first term is a finite summation on which we can apply \autoref{portoy} thus
\[
\sum_{\substack{n\sim n'\\ n\sim q}} \norm{B_{(q,q')}^*B_{(n,n')}}{}{\ha} = O(h^\be).
\]
Let us now consider the second summation. We split it as
\begin{align}\label{3sums}
\sum_{n,n':n\not\sim n'\text{ or }n\not\sim q}\norm{B_{(q,q')}^*B_{(n,n')}}{}{\ha} &= 
\sum_{n,n':\substack{n\not\sim n'\\n\not\sim q}}\norm{B_{(q,q')}^*B_{(n,n')}}{}{\ha} + 
\sum_{n,n':\substack{n\sim n'\\ n\not\sim q}}\norm{B_{(q,q')}^*B_{(n,n')}}{}{\ha} \\&\qquad+ \sum_{n,n':\substack{n\not\sim n'\\ n'\sim q'}}\norm{B_{(q,q')}^*B_{(n,n')}}{}{\ha}.\nonumber
\end{align}
In virtue of \autoref{lemhinty} applied both for $\norm{{}_q B_{+}^*{}_nB_+}{}{}$ and for $\norm{{}_n B_{+}{}_{n'}B_-}{}{}$ for every $k>3$, 
\begin{align*}
\sum_{n,n':\substack{n\not\sim n'\\n\not\sim q}}\norm{B_{(q,q')}^*B_{(n,n')}}{}{\ha} &\le C_kh^{\frac k2(1-\rho)+M(\qt-\frac12\rho)-1} \sum_{n,n':\substack{n\not\sim n'\\n\not\sim q}}\min\left\{\norm{\iota^{-1}(n-q)}{}{-\frac{k}{2}},\norm{\iota^{-1}(n-n')}{}{-\frac{k}{2}}\right\} \\
&\le C_kh^{\frac k2(1-\rho)+M(\qt-\frac12\rho)-1}\cdot \\&\hspace{1em}\int_{\RR^2\setminus \mathbb{B}(q',10\ell(\ka))}\int_{\RR^2\setminus\mathbb{B}(x,10\ell(\ka))}\min\left\{\norm{\iota^{-1}(\xi-q)}{}{-\frac{k}{2}},\norm{\iota^{-1}(\xi`- x)}{}{-\frac{k}{2}}\right\}d\xi dx	 \\ 
&= C_kh^{\frac k2(1-\rho)+M(\qt-\frac12\rho)-1}\iint_{(\RR^2\setminus \mathbb{B}(0,10\ell(\ka)))^2}\min\left\{\norm{\iota^{-1}\xi}{}{-\frac{k}{2}},\norm{\iota^{-1}x}{}{-\frac{k}{2}}\right\}d\xi dx \\
&= O(h^{\frac k2(1-\rho)+M(\qt-\frac12\rho)-1}).
\end{align*}
The other sums in the right hand side of \eqref{3sums} can be treated in a similar manner. First,
\begin{align*}
	\sum_{n,n':\substack{n\sim n'\\ n\not\sim q}}\norm{B_{(q,q')}^*B_{(n,n')}}{}{\ha} &\le C\sum_{n,n':\substack{n\sim n'\\ n\not\sim q}}C_{n,n'}\norm{{}_q B_+ {}_{n}B_+}{}{\ha}\\&\le C_kh^{\frac k2(1-\rho)+M(\qt-\frac12\rho)-1}\sum_{n,n':\substack{n\sim n'\\ n\not\sim q}}C_{n,n'}\norm{\iota^{-1}(n-q)}{}{-\frac{k}{2}}\\&=O(h^{\frac k2(1-\rho)+M(\qt-\frac12\rho)-1}).
\end{align*}
and,
\begin{align*}
\sum_{n,n':\substack{n\not\sim n'\\ n\sim q}}\norm{B_{(q,q')}^*B_{(n,n')}}{}{\ha} &\le C\sum_{n\not\sim n'}\norm{{}_n B_+ {}_{n'}B_-}{}{\ha}\\&\le C_kh^{\frac k2(1-\rho)+M(\qt-\frac12\rho)-1}\sum_{n\not\sim n'}\norm{\iota^{-1}(n-n')}{}{-\frac{k}{2}}\\&=O(h^{\frac k2(1-\rho)+M(\qt-\frac12\rho)-1}).
\end{align*}
The estimates above yield 
\[
\sup_{(q,q')\in\ZZ^4}\sum_{(n,n')\in\ZZ^4}\norm{B_{(q,q')}^*(B_{(n,n')})}{}{\ha} = \sup_{(q,q')\in \ZZ^4}\sum_{(n,n')\in \ZZ^4} \norm{{}_{q'}B^*_-\ {}_{q}B^*_+\ {}_{n}B_+\ {}_{n'}B_-}{}{\ha} = O(h^{\be}).
\]

By analogous methods we obtain as well that 
\[
\sup_{(q,q')\in\ZZ^4}\sum_{(n,n')\in\ZZ^4}\norm{(B_{(q,q')})B_{(n,n')}^*}{}{\ha}=\sup_{(q,q')\in \ZZ^4}\sum_{(n,n')\in \ZZ^4} \norm{{}_{q}B_+\ {}_{q'}B_-\ {}_{n'}B^*_-\ {}_{n}B^*_+}{}{\ha} = O(h^{\be}).
\]
From applying Cotlar-Stein theorem and recalling $\sum_{m\in\ZZ^2}\ind^\ka_{S^m}=1$ $$\sum_{n,n'\in\ZZ^2:\norm{n}{\8}{},\norm{n'}{\8}{}<M} B_{(n,n')}\xrightarrow{M\to\8}\Op_h(b_{w}^+)\Op_h(b_{w}^-)$$ in the strong operator topology and 
\eq{
	\norm{\Op_h(b_{w}^+)\Op_h(b_{w}^-)}{\bltr}{} = O(h^\be).
}{estB}
\end{proof}
\newpage
\bibliographystyle{alpha}
\bibliography{catfup_paper2.bib}{}
\end{document}